\documentclass[reqno]{amsart}
\usepackage[utf8]{inputenc}
\usepackage{amscd,amsfonts,amsmath,amssymb,amsthm,bbm,dsfont,centernot,mathtools,mathrsfs,yhmath,fancyhdr,xy,eso-pic,xcolor,extarrows,relsize,caption,eurosym,tikz,color,mathrsfs,stackrel}
\usepackage[margin=1in]{geometry}
\usepackage[colorlinks=true,linkcolor=blue,citecolor=red]{hyperref}
\usepackage[shortlabels]{enumitem}
\newtheorem{corollary}{Corollary}[section]
\newtheorem{lemma}[corollary]{Lemma}
\newtheorem{proposition}[corollary]{Proposition}
\newtheorem{theorem}[corollary]{Theorem}
\theoremstyle{definition}
\newtheorem{definition}[corollary]{Definition}
\newtheorem{remark}[corollary]{Remark}

\newtheorem*{acknowledgements}{\sc Acknowledgements}
\numberwithin{equation}{section}
\allowdisplaybreaks

\def \div {\mathop {\rm div}\nolimits}
\def \de {\mathrm{d}}
\def \e {\mathrm{e}}
\def \C {\mathbb C}
\def \R {\mathbb R}
\def \N {\mathbb N}
\def \A {\mathbb A}
\def \B {\mathbb B}
\def \O {\Omega}
\def \W {\mathscr{W}}
\def \uez {u^0_{\varepsilon}}
\def \ueu {u^1_{\varepsilon}}
\def \u {u_{\varepsilon}}
\def \du {\dot{u}_{\varepsilon}}
\def \ddu {\ddot{u}_{\varepsilon}}
\def \ube {\Bar{u}_{\varepsilon}}
\def \dube {\dot{\Bar{u}}_{\varepsilon}}

\def \tue {\tilde{u}_{\varepsilon}}
\def \dtue {\dot{\tilde{u}}_{\varepsilon}}

\def \uein {u_{\varepsilon,in}}
\def \duein {\dot{u}_{\varepsilon,in}}
\def \uin {u_{in}}

\def \vez {v^0_{\varepsilon}}
\def \veu {v^1_{\varepsilon}}
\def \v {{{v}}_{\varepsilon}}
\def \dv {\dot{{v}}_{\varepsilon}}
\def \ddv {\ddot{{v}}_{\varepsilon}}
\def \vb {\Bar{v}_{\varepsilon}}
\def \dvb {\dot{\Bar{v}}_{\varepsilon}}

\def \tve {\tilde{v}_{\varepsilon}}
\def \dtve {\dot{\tilde{v}}_{\varepsilon}}
\def \ddtve {\ddot{\tilde{v}}_{\varepsilon}}
\def \he {h_{\varepsilon}}

\def \hh {\hat{h}_{\varepsilon}}
\def \hv {\hat{v}_{\varepsilon}}
\def \hvj {\hat{v}_{\varepsilon_j}}
\def \lee {\ell_{\varepsilon}}
\def \dlee {\dot{\ell}_{\varepsilon}}
\def \gae {\gamma_{\varepsilon}}

\def \fe {f_{\varepsilon}}
\def \phie {\varphi_{\varepsilon}}

\def \gef {g_{\varepsilon}}

\def \tze {\tilde{z}_{\varepsilon}}
\def \zez {z^0_{\varepsilon}}
\def \ze {z_{\varepsilon}}
\def \dze {\dot{z}_{\varepsilon}}
\def \ddze {\ddot{z}_{\varepsilon}}
\def \V {\mathcal{V}}
\def \VD {\mathcal{V}_{0}}
\def \w {{w}_{\varepsilon}}
\def \dw {\dot{{w}}_{\varepsilon}}
\def \wb {\Bar{w}_{\varepsilon}}
\def \dwb {\dot{\Bar{w}}_{\varepsilon}}
\def \uk {u_{\varepsilon_{j}}}
\def \duk {\dot{u}_{\varepsilon_{j}}}
\def \fuk {f_{\varepsilon_{j}}}

\def \sumh {\sigma^{\eta}}
\def \diffh {\delta^{\eta}}

\usepackage{scalerel}
\usepackage[usestackEOL]{stackengine}

\def\dashint{\,\ThisStyle{\ensurestackMath{%
            \stackinset{c}{.2\LMpt}{c}{.5\LMpt}{\SavedStyle-}{\SavedStyle\phantom{\int}}}%
        \setbox0=\hbox{$\SavedStyle\int\,$}\kern-\wd0}\int}

\begin{document}

\title{Quasistatic limit of a dynamic viscoelastic model with memory}
\author[G. Dal Maso]{Gianni Dal Maso}
\address[Gianni Dal Maso]{SISSA, via Bonomea 265, 34136 Trieste, Italy}
\email{dalmaso@sissa.it}
\author[F. Sapio]{Francesco Sapio}
\address[Francesco Sapio]{SISSA, via Bonomea 265, 34136 Trieste, Italy}
\email{fsapio@sissa.it}

\begin{abstract}
We study the behaviour of the solutions to a dynamic evolution problem for a viscoelastic model with long memory, when the rate of change of the data tends to zero. We prove that a suitably rescaled version of the solutions converges to the solution of the corresponding stationary problem.
\end{abstract}

\maketitle


\medskip


\section{Introduction}

The most common models of viscoelasticity with long memory, such as the Maxwell model (see \cite{DL}, \cite{Fab}, \cite{MF-AM}, \cite{SL}), lead to a dynamic evolution governed by a system of partial differential equations of the form
\begin{equation}\label{DM-viscoel-dyn}
    \ddot u(t)-\div((\A+\B)eu(t))+\int_{-\infty}^t\frac{1}{\beta}\e^{-\frac{t-\tau}{\beta}}\div(\mathbb B e u(\tau))\de \tau=\ell(t)\quad\text{in $\O$ for $t\in[0,T]$},
\end{equation}
where $\Omega\subset\R^d$ is the reference configuration, $[0,T]$ is the time interval, $u(t)$ and $eu(t)$ are the displacement at time $t$ and the symmetric part of its gradient, $\A$ and $\B$ are the elasticity and viscosity tensors, $\beta>0$ is a material constant, and
$\ell(t)$ is the external load at time $t$. 
This system is complemented by boundary and initial conditions
\begin{alignat}{3}
&u(t)=z(t) &&  \quad\text{on $\partial\O$ for $t\in[0,T]$},\label{DM-bound1}\\
&u(t)=\uin(t)&&\quad \text{in $\O$ for $t\in(-\infty,0]$},\label{DM-initial}
\end{alignat}
where $z$ and $\uin$ are prescribed functions, the latter representing the history of the displacement for $t\leq 0$. Existence and uniqueness for \eqref{DM-viscoel-dyn}--\eqref{DM-initial} can be found in \cite{DA}.

In this paper we study the quasistatic limit of the solutions to this problem, i.e., the limit of these solutions when the rate
of change of the data tends to zero.
More precisely, given a small parameter $\varepsilon>0$, we consider the solution 
$u^{\varepsilon}$ of \eqref{DM-viscoel-dyn}--\eqref{DM-initial} corresponding to
$\ell(\varepsilon t)$, $z(\varepsilon t)$, and $\uin(\varepsilon t)$.
To study the asymptotic behaviour of $u^{\varepsilon}$ as $\varepsilon\to0^+$ it is convenient to
introduce the rescaled solution $\u(t):=u^{\varepsilon}(t/\varepsilon)$, which
turns out to be the solution of the system
\begin{equation}\label{DM-viscoel-dyn-e}
    \varepsilon^2\ddu(t)-\div((\A+\B)e\u(t))+\int_{-\infty}^t\frac{1}{\beta\varepsilon}\e^{-\frac{t-\tau}{\beta\varepsilon}}\div(\mathbb B e \u(\tau))\de \tau=\ell(t)\quad\text{in $\O$ for $t\in[0,T]$},
\end{equation}
with boundary and initial conditions \eqref{DM-bound1} and \eqref{DM-initial}.

Under different assumptions on $\ell(t)$, $z(t)$, and $\uin(t)$ we prove (Theorems \ref{main-thm-f+g} and \ref{main-f}) that $\u(t)$ converges, as $\varepsilon\to 0^+$, to the solution $u_0(t)$ of the stationary problem
\begin{equation}\label{DM-elliptic}
    -\div(\A eu_0(t))=\ell(t)\quad\text{in $\O$ for $t\in[0,T]$}, 
\end{equation}
with boundary condition \eqref{DM-bound1}. 

It is not difficult to prove a similar result for the Kelvin-Voigt model, in which the viscosity term
\begin{equation}\label{DM-damping}
-\div(\B eu(t))+\int_{-\infty}^t\frac{1}{\beta}\e^{-\frac{t-\tau}{\beta}}\div(\mathbb B e u(\tau))\de \tau
\end{equation}
is replaced by $-\div(\B e\dot{u}(t))$. On the other hand, it is well known that, in general, the convergence of $\u$ to $u_0$ does not hold for the equation of elastodynamics without damping terms, i.e., when $\B=0$. The purpose of this paper is to prove that
the non-local damping term \eqref{DM-damping} is enough to obtain the convergence of the solutions of the evolution problems to the solution of the stationary problem. 

The main tools to prove our results are two different estimates (Lemmas \ref{lem:en-est} and \ref{lem:ex-lax}), related to the energy-dissipation balance \eqref{energy} and to the elliptic system \eqref{weak-fomulation-lap} obtained from \eqref{DM-viscoel-dyn-e} via Laplace Transform. After a precise statement of all assumptions, more details on the line of proof will be given after Theorem~\ref{main-f}.

\section{Hypotheses and statement of the problem}

Let $d$ be a positive integer and let $\O\subset\R^d$ be a bounded open set with Lipschitz boundary.
We use standard notation for Lebesgue and Sobolev spaces. Let $\mathbb R^{d\times d}_{sym}$ be the space of all symmetric $d{\times}d$ matrices. For convenience we set
\begin{equation}\label{spaziQ}
H:=L^2(\O;\R^d),\quad \tilde{H}:=L^2(\O;\R^{d\times d}_{sym}),\quad V:=H^1(\O;\R^d),\quad V_0:=H_0^1(\O;\R^d),\quad V'_0:=H^{-1}(\O;\R^d),
\end{equation}
and we always identify the dual of $H$ with $H$ itself. The symbols $(\cdot,\cdot)$ and $\|\cdot\|$ denote the scalar product and the norm in $H$ or in $\tilde{H}$, according to the context. The duality product between $V'_0$ and $V_0$ is denoted by $\langle\cdot,\cdot\rangle$.
Given $u\in V$, its strain $eu$ is defined as the symmetric part of its gradient, i.e., $eu:=\frac{1}{2}(\nabla u+\nabla u^T)$, where $\nabla u$ is the Jacobian matrix, whose
components are
$(\nabla u)_{ij}:= \partial_j u_i$ for $i,j=1,\dots,d$.

Under these assumptions, the Second Korn Inequality (see, e.g.,~\cite[Theorem~2.4]{OSY}) states that there exists a positive constant $C_K=C_K(\O)$ such that 
\begin{equation}\label{eq:hyp_kornQ}
\|\nabla u\|\leq C_K\left(\|u\|^2+\|e u\|^2\right)^{1/2}\quad\text{for every }u\in V.
\end{equation} 
Moreover, there exists a positive constant $C_{P}=C_P(\O)$ such that the following Korn-Poincaré Inequality holds (see, e.g.,~\cite[Theorem~2.7]{OSY}):  
\begin{equation}\label{KP-ineq}
\|u\|\leq C_{P}\|e u\|\qquad\text{for every }u\in V_0.
\end{equation}
Thanks to \eqref{eq:hyp_kornQ} we can use on the space $V$ the equivalent norm 
\begin{equation*}
\|u\|_{V}:=(\|u\|^2+\|e u\|^2)^{1/2}\quad\text{for every }u\in V.
\end{equation*}

Let $\mathscr L(\R^{d\times d}_{sym};\R^{d\times d}_{sym})$ be the space of all linear operators from $\R^{d\times d}_{sym}$ into itself. We assume that the elasticity and viscosity tensors $\A$ and $\B$, which depend on the variable $x\in\O$, satisfy the following assumptions:
\begin{align}
    &\A,\B\in L^{\infty}(\O;\mathscr L(\R^{d\times d}_{sym};\R^{d\times d}_{sym})),\label{CB1Q}\\
    &\A(x)\xi_1\cdot\xi_2=\xi_1\cdot \A(x)\xi_2,&&\B(x)\xi_1\cdot\xi_2=\xi_1\cdot \B(x)\xi_2&&\hspace{-2pt}\text{for a.e.\ $x\in\O$ and every $\xi_1,\xi_2\in\R_{sym}^{d\times d}$},\label{sym}\\
    &c_{\A}|\xi|^2\leq\A(x)\xi\cdot\xi\leq C_{\A}|\xi|^2,&&c_{\B}|\xi|^2\leq\B(x)\xi\cdot\xi\leq C_{\B}|\xi|^2&&\hspace{-2pt}\text{for a.e.\ $x\in\O$ and every $\xi\in\R_{sym}^{d\times d}$,}\label{CB2Q}
\end{align}
where $c_{\A}$, $c_{\B}$, $C_{\A}$, and $C_{\B}$ are positive constants independent of $x$, and the dot denotes the Euclidean scalar product of matrices.

Let us fix $T>0$ and $\beta>0$. To give a precise meaning to the notion of solution to problem \eqref{DM-bound1}--\eqref{DM-viscoel-dyn-e} we introduce the function spaces
\begin{align*}
\V:=L^2(0,T;V)\hspace{2pt}\cap &\hspace{2pt}H^1(0,T;H)\cap H^2(0,T;V'_{0}),\quad \VD:=L^2(0,T;V_0)\cap H^1(0,T;H)\cap H^2(0,T;V'_{0}),
\\
&\V_{loc}:=L^2_{loc}(-\infty,T;V)\cap H^1_{loc}(-\infty,T;H)\cap H^2_{loc}(-\infty,T;V'_{0}).
\end{align*}

\begin{remark}
By the Sobolev Embedding Theorem,
if $u\in\V$ (resp.\ $u\in\V_{loc}$), then 
$$u\in C^0([0,T];H)\cap C^1([0,T];V_0')\quad\text{(resp.\ $u\in C^0((-\infty,T);H)\cap C^1((-\infty,T);V_0')$)}.$$
\end{remark}

We study problem \eqref{DM-bound1}--\eqref{DM-viscoel-dyn-e} with $\ell$, $z$, and $\uin$ depending on $\varepsilon$. Let us consider $\varepsilon>0$ and 
\begin{equation}\label{hyp1}
    \fe\in L^2(0,T;H),\quad \gef\in H^1(0,T;V'_0),\quad \ze \in H^2(0,T;H)\cap H^1(0,T;V),
\end{equation}
$\uein\in C^0((-\infty,T);H)\cap C^1((-\infty,T);V_0')$ such that
\begin{align}
    \uein(0)\in V,&\quad \uein(0)-\ze(0)\in V_0, \quad \duein(0)\in H,\quad \int_{-\infty}^0\frac{1}{\beta\varepsilon}\e^{\frac{\tau}{\beta\varepsilon}}\|\uein(\tau)\|_V\de\tau< +\infty\label{hyp3}.
\end{align}

The notion of solution to \eqref{DM-bound1}--\eqref{DM-viscoel-dyn-e} is made precise by the following definition.

\begin{definition}\label{def6666}
We say that $\u$ is a solution to the viscoelastic dynamic system \eqref{DM-bound1}--\eqref{DM-viscoel-dyn-e}, with forcing term $\ell=\fe+\gef$, boundary condition $\ze$, and initial condition $\uein$, if
\begin{equation}\label{spaz}
    \hspace{-0.1cm}\begin{cases}
    \u\in\V_{loc}\quad\text{and}\quad\u-\ze\in\VD,\\
    \displaystyle\varepsilon^2\ddu(t) -\div((\A+\B) e\u(t))+\int_{-\infty}^t\frac{1}{\beta\varepsilon}\e^{-\frac{t-\tau}{\beta\varepsilon}}\div(\B e\u(\tau))\de\tau=\fe(t)+\gef(t)\quad\hspace{-0.1cm}\text{for a.e.\ }t\in [0,T], \\
\u(t)=\uein(t)\quad\text{for every }t\in(-\infty,0].
    \end{cases}
\end{equation}
\end{definition}

In the next remark we shall see that \eqref{spaz} can be reduced to the following problem starting from $0$:
\begin{equation}\label{spaz2}
    \begin{cases}
    \u\in\V\quad\text{and}\quad\u-\ze\in\VD,\\
     \displaystyle\varepsilon^2\ddu(t) -\div((\A+\B) e\u(t))+\int_0^t\frac{1}{\beta\varepsilon}\e^{-\frac{t-\tau}{\beta\varepsilon}}\div(\B e\u(\tau))\de\tau= \phie(t)+ \gae(t)\quad\hspace{-0.1cm}\text{for a.e.\ }t\in [0,T],\\
     \u(0)=\uez\text{ in }H \quad \text{and}\quad \du(0)=\ueu\text{ in }V_0',
    \end{cases}
\end{equation}
with $\phie\in L^2(0,T;H)$, $\gae\in H^1(0,T;V'_0)$, $\uez\in V$, $\uez-\ze(0)\in V_0$, $\ueu\in H$.

\begin{remark}\label{corris}
It is easy to see that $\u$ is a solution according to Definition~\ref{def6666} if and only if its restriction  to $[0,T]$, still denoted by $\u$, solves \eqref{spaz2} with
\begin{equation}
    \phie=\fe, \quad\gae=\gef-p_{\varepsilon},\quad \uez=\uein(0),\quad\ueu=\duein(0),\label{u1}
\end{equation}
where
\begin{equation}\label{g0}
    p_{\varepsilon}(t):=\e^{-\frac{t}{\beta\varepsilon}}g^0_{\varepsilon}\quad\text{with}\quad g^0_{\varepsilon}:=\int_{-\infty}^0\frac{1}{\beta\varepsilon}\e^{\frac{\tau}{\beta\varepsilon}}\div(\B e\uein(\tau))\de\tau.
\end{equation}
\end{remark}

To solve problem \eqref{spaz2} it is enough to study the corresponding problem with 
homogeneous boundary condition:
\begin{equation}\label{weaksol2}
\begin{cases}
\v\in \VD,
\\
\displaystyle
\varepsilon^2\ddv(t)- \div((\A+\B) e\v(t))+\int_0^t\frac{1}{\beta\varepsilon}\e^{-\frac{t-\tau}{\beta\varepsilon}}\div(\B e\v(\tau))\de\tau=\he(t)+\lee(t)\quad\text{for a.e.\ }t\in [0,T],
\\
\v(0)= \vez\text{ in }H \quad \text{and}\quad \dv(0)=\veu\text{ in }V_0',
\end{cases}
\end{equation}
with 
\begin{equation}\label{hyprr}
    \he\in L^2(0,T;H),\quad \lee\in H^1(0,T;V'_0),\quad \vez\in V_0,\quad \veu\in H.
\end{equation}
\begin{remark}\label{rem-um}
The function $\u$ is a solution to \eqref{spaz2} if and only if $\v=\u-\ze$ solves \eqref{weaksol2}
with
\begin{align}
\he(t)=\phie(t)-\varepsilon^2\ddze(t),\quad &\lee(t)=\gae(t)+\div((\A+\B)e\ze(t))-\int_0^t\frac{1}{\beta\varepsilon}\e^{-\frac{t-\tau}{\beta\varepsilon}}\div(\B e\ze(\tau))\de\tau,\nonumber\\
&\vez=\uez-\ze(0),\quad \veu=\ueu-\dze(0),\label{datin}
\end{align}
Therefore, existence and uniqueness for \eqref{weaksol2} imply existence and uniqueness for \eqref{spaz2}.
\end{remark}

\begin{remark}\label{rem-un-2}
In \cite{DA} problem \eqref{weaksol2} has been studied with initial conditions taken in the sense of interpolation spaces. Given two Hilbert spaces $X$ and $Y$, the symbol $[X,Y]_{\theta}$ denotes the interpolation space between $X$ and $Y$ of exponent $\theta\in(0,1)$. Thanks to \cite[Theorem 3.1]{LM} we have the following inclusions:
\begin{equation*}
    L^2(0,T;V_0)\cap H^1(0,T;H)\subset C^0([0,T];V_0^{\frac{1}{2}})\qquad\text{and}\qquad L^2(0,T;H)\cap H^1(0,T;V'_0)\subset C^0([0,T];V_0^{-\frac{1}{2}}),
\end{equation*}
where $V_0^{\frac{1}{2}}:=[V_0,H]_{\frac{1}{2}}$ and $V_0^{-\frac{1}{2}}:=[H,V'_0]_{\frac{1}{2}}$.  Consequently
\begin{equation*}
    \VD\subset  C^0([0,T];V_0^{\frac{1}{2}})\cap C^1([0,T];V_0^{-\frac{1}{2}}).
\end{equation*}
Therefore, the initial conditions in \eqref{weaksol2} are satisfied also in the stronger sense
\begin{equation}\label{dat-um}
\v(0)=\vez\text{ in }V_0^{\frac{1}{2}} \quad \text{and}\quad \dv(0)=\veu\text{ in }V_0^{-{\frac{1}{2}}}.
\end{equation}
\end{remark}

The following proposition provides the main properties of the solutions. 
We recall that, if $X$ is a Banach space, $C_w^0([0,T];X)$ denotes the space of all weakly continuous functions from $[0,T]$ to $X$, namely, the vector space of all functions $u\colon [0,T]\to X$ such that for every $x'\in X'$ the function $t\mapsto \langle x', u(t)\rangle$ is continuous from $[0,T]$ to $\R$.

\begin{proposition}\label{en-eq}
Given $\varepsilon>0$, assume \eqref{hyp1} and \eqref{hyp3}. Then there exists a unique solution $\u$ to the viscoelastic dynamic system \eqref{spaz}. Moreover, it satisfies
\begin{equation}\label{cont}
    \u\in C^0([0,T];V)\cap C^1([0,T];H).
\end{equation}
\end{proposition}

\begin{proof}
By Remarks \ref{corris} and \ref{rem-um} it is enough to prove the theorem for \eqref{weaksol2}. Existence and uniqueness are proved in \cite{DA}, taking into account Remark \ref{rem-un-2} about the equivalence between the initial conditions in the sense of \eqref{weaksol2} and \eqref{dat-um}. 

After an integration by parts with respect to time, it easy to see that the weak formulation \eqref{weaksol2} is equivalent to the following one: 
\begin{align}\label{weak-weak-form-F}
-\varepsilon^2\int_0^T(\dv(t),\dot \varphi(t))\de t&+\int_0^T((\A+\B) e\v(t),e \varphi(t))\de t-\int_0^T\int_0^t\frac{1}{\beta\varepsilon}\e^{-\frac{t-\tau}{\beta\varepsilon}}(\B e\v(\tau),e \varphi(t))\de\tau\de t\nonumber\\
&=\int_0^T( \he(t),\varphi(t))\de t+\int_0^T\langle \lee(t),\varphi(t)\rangle\de t\quad \text{for every $\varphi \in C^{\infty}_c(0,T;V)$.}
\end{align}
In \cite{S}, in a more general context, it has been proved that if $\v$ satisfies \eqref{weak-weak-form-F} and the initial conditions in the sense of \eqref{weaksol2}, then it satisfies also 
\begin{align}
     \v\in C_w^0([0,T];V) \quad&\text{and}\quad \dv\in C_w^0([0,T];H),\label{cont-weak}\\
    \lim_{t\to 0^+}\|\v(t)-\vez\|_V=0\quad &\text{and}\quad \lim_{t\to 0^+}\|\dv(t)-\veu\|=0\nonumber.
\end{align}

We fix $s\in [0,T)$. We want to prove 
\begin{equation}\label{cont-right}
    \lim_{t\to s^+}\|\v(t)-\v(s)\|_V=0\qquad\text{and}\qquad  \lim_{t\to s^+}\|\dv(t)-\dv(s)\|=0.
\end{equation}
Thanks to the theory developed in \cite{DA} there exists a unique $\tve\in L^2(s,T;V_0)\cap H^1(s,T;H)\cap H^2(s,T;V'_{0})$ such that
\begin{align}
\varepsilon^2\ddtve(t) &-\div((\A+\B) e\tve(t))+\int_s^t\frac{1}{\beta\varepsilon}\e^{-\frac{t-\tau}{\beta\varepsilon}}\div(\B e\tve(\tau)) \de\tau\nonumber\\
&=\he(t)+\lee(t)-\int_0^s\frac{1}{\beta\varepsilon}\e^{-\frac{t-\tau}{\beta\varepsilon}}\div(\B e\v(\tau))\de\tau\quad\text{for a.e.\ $t\in [s,T]$},\label{eqnew}
\end{align}
\begin{equation}\label{eqnew2}
    \lim_{t\to s^+}\|\tve(t)-\v(s)\|=0\quad \text{and}\quad \lim_{t\to s^+}\|\dtve(t)-\dv(s)\|_{V'_0}=0.
\end{equation}
By the results in \cite{S} the function $\tve$ satisfies also
\begin{equation}\label{bis}
\lim_{t\to s^+}\|\tve(t)-\v(s)\|_V=0\quad \text{and}\quad \lim_{t\to s^+}\|\dtve(t)-\dv(s)\|=0.
\end{equation}
Since clearly $\v$ satisfies \eqref{eqnew} and \eqref{eqnew2}, by uniqueness we have $\tve(t)=\v(t)$ for every $t\in[s,T]$. In particular, from \eqref{bis} we deduce that \eqref{cont-right} holds.
\phantom\qedhere
\end{proof}

To complete the proof we need the following proposition about the energy-dissipation balance, where $\tilde{H}$ is defined by \eqref{spaziQ}, and $\W_{\varepsilon}(t)$ represents the work done in the interval $[0,t]$.

\begin{proposition}\label{diss-bal}
Given $\varepsilon>0$, we assume \eqref{hyprr}. Let $\v$ be the solution to \eqref{weaksol2} and let $\w\colon [0,T]\rightarrow \tilde{H}$ be defined by 
\begin{equation}\label{we}
    \w(t):=\e^{-\frac{t}{\beta\varepsilon}}\int_0^t\frac{1}{\beta\varepsilon}\e^{\frac{\tau}{\beta\varepsilon}}e\v(\tau)\de \tau\quad\text{for every $t\in[0,T]$.}
\end{equation}
Then $\w\in H^1(0,T;\tilde{H})$ and the following energy-dissipation balance holds for every $t\in[0,T]$:
\begin{align}\label{energy}
    \frac{\varepsilon^2}{2}\|\dv(t)\|^2&+\frac{1}{2}(\A e\v(t),e\v(t))+\frac{1}{2}(\B (e\v(t)-\w(t)),e\v(t)-\w(t))\nonumber\\
    &+\beta\varepsilon\int_0^t(\B \dw(\tau),\dw(\tau))\de \tau= \frac{\varepsilon^2}{2}\|\veu\|^2+\frac{1}{2}((\A+\B) e\vez,e\vez)+\W_{\varepsilon}(t),
\end{align}
where 
\begin{align*}
  \W_{\varepsilon}(t):&=\int_0^t(\he(\tau),\dv(\tau))\de \tau-\int_0^t \langle \dlee(\tau),\v(\tau)\rangle\de\tau+\langle \lee(t),\v(t)\rangle-\langle \lee(0),\vez\rangle.
  \end{align*}
\end{proposition}

\begin{proof}
It is convenient to extend the data of our problem to the interval $[0,2T]$ by setting
\begin{equation*}
    \he(t):=0\quad\text{and}\quad \lee(t):=\lee(T)\quad \text{for every $t\in (T,2T]$}.
\end{equation*}
It is clear that $\he\in L^2(0,2T;H)$ and $\lee\in H^1(0,2T;V'_0)$. By uniqueness of the solution to \eqref{weaksol2}, the solution on $[0,2T]$ is an extension of $\v$, still denoted by $\v$. We also consider the extension of $\w$ on $[0,2T]$ defined by \eqref{we}. 

Since $e\v\in L^2(0,2T;\tilde H)$, it follows from \eqref{we} that $\w\in H^1(0,2T;\tilde{H})$, and
\begin{equation}\label{ord}
    \beta\varepsilon \dw(t)=e\v(t)-\w(t)\quad \text{for a.e.\ $t\in[0,2T]$}.
\end{equation}
Thanks to \eqref{cont-right} in $[0,2T]$ and \eqref{ord} there exists a representative of $\dw$ such that
\begin{equation}\label{cont-dw}
    \lim_{t\to s^+}\|\dw(t)-\dw(s)\|=0\quad\text{for every $s\in[0,2T)$.}
\end{equation}
Moreover, since $\v$ satisfies \eqref{weaksol2} in $[0,2T]$, we have
\begin{equation}\label{equiv-eq}
\varepsilon^2\ddv(t)-\div(\A e\v(t))-\div(\B (e\v(t)-\w(t)))= \he(t)+\lee(t)\quad\text{for a.e.\ $t\in[0,2T]$}.
\end{equation}
Multiplying \eqref{ord} and \eqref{equiv-eq} by $\psi\in \tilde{H}$ and $\varphi\in V_0$, respectively, and then integrating over $\O$ and adding the results, for a.e.\ $t\in[0,2T]$ we get 
\begin{align}\label{eq-finale}
\varepsilon^2\langle\ddv(t), \varphi\rangle +(\A e\v(t),e \varphi)+(\B (e\v(t)&-\w(t)),e \varphi-\psi)+\beta\varepsilon(\B \dw(t),\psi)=( \he(t),\varphi)+\langle \lee(t),\varphi\rangle.
\end{align}

Given a function $r$ from $[0, 2T]$ into a Banach space $X$, for every $\eta > 0$ we define the sum and the difference $\sumh r, \diffh r\colon [0, 2T -\eta]\rightarrow X $ by $\sumh r(t):= r(t+\eta)+r(t)$ and $\diffh r(t):= r(t+\eta)-r(t)$. For a.e.\ $t\in [0,2T-\eta]$ we have $\sumh \v(t), \diffh \v(t)\in V_0$ and $\sumh \w (t), \diffh\w (t)\in \tilde{H}$. For a.e.\ $t\in[0,2T-\eta]$ we use \eqref{eq-finale} first at time $t$ and then at time $t + \eta$, with $\varphi:=\diffh \v(t)$ and $\psi:=\diffh \w (t)$. By summing the two expressions and then integrating in time on the interval $[0,t]$ we get
\begin{equation}\label{t+h,h}
  \int_0^t[\varepsilon^2 K_{\eta}(\tau)+A_{\eta}(\tau)+B_{\eta}(\tau)+\varepsilon D_{\eta}(\tau)]\de\tau=\int_0^t W_{\eta}(\tau)\de\tau,
\end{equation}
where for a.e.\ $\tau\in[0,2T-\eta]$ 
\begin{align}
    & K_{\eta}(\tau):=\langle\sumh \ddv(\tau), \diffh  \v(\tau)\rangle,\nonumber\\
    & A_{\eta}(\tau):=(\A\, \sumh e\v(\tau),\diffh  e\v(\tau)),\nonumber\\
    & B_{\eta}(\tau):=(\B (\sumh e\v(\tau)- \sumh \w(\tau)),\diffh  e\v(\tau)-\diffh  \w(\tau)),\nonumber\\
    & D_{\eta}(\tau):=\beta(\B\, \sumh\dw(\tau),\diffh  \w (\tau)),\nonumber\\
    & W_{\eta}(\tau):=( \sumh\he(\tau),\diffh  \v(\tau))+\langle\sumh \lee(\tau),\diffh  \v(\tau)\rangle.\nonumber
\end{align}
An integration by parts in time gives
\begin{align}\label{int-K}
    \int_0^tK_{\eta}(\tau)\de \tau&=(\sumh\dv(t), \diffh  \v(t))-(\sumh\dv(0), \diffh  \v(0))-\int_0^t(\sumh\dv(\tau), \diffh  \dv(\tau))\de \tau\nonumber\\
   &=\int_t^{t+\eta}(\sumh\dv(t), \dv(\tau))\de \tau-\int_0^{\eta}(\sumh\dv(0),\dv(\tau) )\de \tau-\int_0^{t}\|\dv(\tau+h)\|^2\de\tau+\int_0^{t}\|\dv(\tau)\|^2\de\tau\nonumber\\
   &=\int_t^{t+\eta}\big[(\sumh\dv(t), \dv(\tau))-\|\dv(\tau)\|^2\big]\de \tau-\int_0^{\eta}\big[(\sumh\dv(0),\dv(\tau) )-\|\dv(\tau)\|^2\big]\de\tau.
\end{align}
Moreover
\begin{align}
&\hspace{-5pt}\int_0^t\hspace{-2pt} A_{\eta}(\tau)\de \tau=\int_t^{t+{\eta}}(\A e\v(\tau),e\v(\tau))\de \tau-\int_0^{{\eta}}(\A e\v(\tau),e\v(\tau))\de \tau,\\
&\hspace{-5pt}\int_0^t \hspace{-2pt}B_{\eta}(\tau)\de \tau=\int_t^{t+{\eta}}\hspace{-2pt}(\B( e\v(\tau)\hspace{-1.5pt}-\hspace{-1.5pt}\w(\tau)),e\v(\tau)\hspace{-1.5pt}-\hspace{-1.5pt}\w(\tau))\de \tau\hspace{-0.5pt}-\hspace{-0.5pt}\int_0^{{\eta}}\hspace{-1pt}(\B( e\v(\tau)\hspace{-1.5pt}-\hspace{-1.5pt}\w(\tau)),e\v(\tau)\hspace{-1.5pt}-\hspace{-1.5pt}\w(\tau))\de \tau,\\
 &\hspace{-5pt} \int_0^t\hspace{-2pt} D_{\eta}(\tau)\de\tau=\beta\int_0^t\int_{\tau}^{\tau+{\eta}}(\B\, \sumh\dw(\tau),\dw (s))\de s\de \tau,\\
 &\hspace{-5pt}\int_0^t\hspace{-2pt} W_{\eta}(\tau)\de\tau=\int_0^t\int_{\tau}^{\tau+{\eta}}( \sumh\he(\tau),\dv(s))\de s\de\tau-\int_{\eta}^t\int_{\tau-{\eta}}^{\tau+{\eta}}\langle \dlee(s),\v(\tau)\rangle\de s\de\tau\nonumber\\
&\hspace{3cm}+\int_{t-{\eta}}^{t}\langle \sumh\lee(\tau),\v(\tau+{\eta})\rangle\de\tau-\int_0^{{\eta}}\langle \sumh\lee(\tau),\v(\tau)\rangle\de\tau.\label{work}
\end{align}
We now divide by ${\eta}$ all terms of \eqref{int-K}--\eqref{work}. Observing that
\begin{align*}
&\sumh \he\xrightarrow[{\eta}\to 0^+]{}2\he\quad \text{strongly in $L^2(0,T;H)$,}\\
&\int_0^t\Big\|\dashint_{\tau}^{\tau+{\eta}}\dv(s)\de s-\dv(\tau)\Big\|^2\de\tau\xrightarrow[{\eta}\to 0^+]{}0,\\
&\int_{\eta}^t\Big\|\dashint_{\tau-{\eta}}^{\tau+{\eta}}\dlee(s)\de s-\dlee(\tau)\Big\|_{V'_0}^2\de\tau\xrightarrow[{\eta}\to 0^+]{}0,
\end{align*}
thanks to \eqref{cont-right} in $[0,2T)$ and \eqref{cont-dw}, we can pass to the limit as ${\eta}\to 0^+$, and from \eqref{t+h,h} we obtain that \eqref{energy} is satisfied for every $t\in [0,T]$. 
\end{proof}

\begin{proof}[Proof of Proposition \ref{en-eq} (Continuation)]
Now we want to prove \eqref{cont}. By using \eqref{energy}, for every $t\in[0,T]$ we can write
\begin{align}\label{newenergy}
\frac{\varepsilon^2}{2}\|\dv(t)\|^2+\frac{1}{2}((\A &+\B)  e\v(t),e\v(t))=\frac{\varepsilon^2}{2}\|\veu\|^2+\frac{1}{2}((\A +\B)  e\vez,e\vez)+ \W_{\varepsilon}(t)\nonumber\\
&-\frac{1}{2}(\B \w(t),\w(t))+(\B e\v(t),\w(t))-\beta\varepsilon\int_0^t(\B\dw(\tau),\dw(\tau))\de\tau.
\end{align}
Let $\Psi_{\varepsilon}\colon[0,T]\rightarrow [0,+\infty)$ be defined by
\begin{equation*}
    \Psi_{\varepsilon}(t):=\frac{\varepsilon^2}{2}\|\dv(t)\|^2+\frac{1}{2}((\A +\B)  e\v(t),e\v(t));
\end{equation*}
since $\w\in C^0([0,T];\tilde{H})$, thanks to \eqref{cont-weak} and \eqref{newenergy} we have $\Psi_{\varepsilon}\in C^0([0,T])$. 

Now we fix $t\in[0,T]$. Given a sequence $\{t_k\}_k\subset [0,T]$ such that $t_k\rightarrow t$ as $k\to +\infty$, we define
\begin{equation*}
    \mathscr{E}_k:=\frac{\varepsilon^2}{2}\|\dv(t_k)-\dv(t)\|^2+\frac{1}{2}((\A +\B)  (e\v(t_k)-e\v(t)),e\v(t_k)-e\v(t)).
\end{equation*}
By elementary computations we have
\begin{equation*}
    \mathscr{E}_k=\Psi_{\varepsilon}(t_k)+\Psi_{\varepsilon}(t)-\varepsilon^2(\dv(t_k),\dv(t))-((\A +\B)  e\v(t_k),e\v(t)),
\end{equation*}
therefore, by \eqref{KP-ineq} and \eqref{CB2Q} there exists a positive constant $C=C(\A,\B,\O)$ such that
\begin{align*}
    \varepsilon^2\|\dv(t_k)-\dv(t)\|^2&+\|\v(t_k)-\v(t)\|^2_V\nonumber\\
    &\leq C\Big(\Psi_{\varepsilon}(t_k)+\Psi_{\varepsilon}(t)-\varepsilon^2(\dv(t_k),\dv(t))-((\A +\B)  e\v(t_k),e\v(t))\Big).
\end{align*}
The right-hand side of the previous inequality tends to $0$ as $k\to +\infty$ because of \eqref{cont-weak} and the continuity of $\Psi_{\varepsilon}$. Since $\ze\in C^0([0,T];V)$, by \eqref{hyp1}, and $\u=\v+\ze$, we obtain \eqref{cont}.
\end{proof}

\section{Statement of the main results}

In this section we present the main results about the convergence, as $\varepsilon\to 0^+$, of the solutions $\u$. We assume the following hypotheses on the dependence on $\varepsilon>0$ of our data:
\begin{itemize}
    \item[(H1)] $\{\fe\}_{\varepsilon}\subset L^2(0,T;H)$, $f\in L^2(0,T;H)$, $\{\gef\}_{\varepsilon}\subset H^1(0,T;V'_0)$, $g\in W^{1,1}(0,T;V'_0)$,     \begin{equation*}
        \fe\xrightarrow[\varepsilon\to 0^+]{}f\quad\quad\text{strongly in $L^2(0,T;H)$},\quad\text{and}\quad \quad\gef\xrightarrow[\varepsilon\to 0^+]{}g\quad\text{strongly in $W^{1,1}(0,T;V'_0)$};
    \end{equation*}
    \item[(H2)] $\{\ze\}_{\varepsilon}\subset H^2(0,T;H)\cap H^1(0,T;V)$, $z\in W^{2,1}(0,T;H)\cap W^{1,1}(0,T;V)$, and
    \begin{equation*}
        \ze\xrightarrow[\varepsilon\to 0^+]{}z\quad\text{strongly in $W^{2,1}(0,T;H)\cap W^{1,1}(0,T;V)$;}
        \end{equation*}
    \item[(H3)] $\{\uein\}_{\varepsilon}\subset C^0((-\infty,0];V)\cap C^1((-\infty,0];H)$, $\uin\in C^0((-\infty,0];V)$, and there exist $a>0$ such that
    \begin{align*}
        &\uein\xrightarrow[\varepsilon\to 0^+]{}\uin\quad\text{strongly in $C^{0}([-a,0];V)$},&&\qquad \varepsilon\duein\xrightarrow[\varepsilon\to 0^+]{}0\quad\text{strongly in $C^{0}([-a,0];H)$},\\
        &\int_{-\infty}^{-a}\frac{1}{\beta\varepsilon}\e^{\frac{\tau}{\beta\varepsilon}}\|\uein(\tau)\|_V\de\tau\xrightarrow[\varepsilon\to 0^+]{}0, &&\qquad \int_{-\infty}^{-a}\frac{1}{\beta\varepsilon}\e^{\frac{\tau}{\beta\varepsilon}}\|\uin(\tau)\|_V\de\tau\xrightarrow[\varepsilon\to 0^+]{}0.
    \end{align*}
\end{itemize} 

\begin{remark}
Let $\uez=\uein(0)$, $\ueu=\duein(0)$, and $u^0=\uin(0)$. Hypothesis (H3) implies 
\begin{equation*}
    \uez\xrightarrow[\varepsilon\to 0^+]{}u^0\quad\text{strongly in $V$}\quad\text{and}\quad  \varepsilon\ueu\xrightarrow[\varepsilon\to 0^+]{}0\quad\text{strongly in $H$}.
\end{equation*}
\end{remark}

Our purpose is to show that the solutions $\u$ converge, as $\varepsilon\to 0^+$, to the solution $u_0$ of the stationary problem \eqref{DM-elliptic} with 
boundary condition \eqref{DM-bound1}. The notion of solution to this problem is the usual one:
\begin{equation}\label{def-ellip}
\begin{cases}
u_0(t)\in V,\quad u_0(t)-z(t)\in V_0,&\quad \text{for a.e. $t\in [0,T]$},\\
-\div(\A eu_0(t))=f(t)+g (t)&\quad \text{for a.e. $t\in [0,T]$}.
\end{cases}
\end{equation}

\begin{remark}
The existence and uniqueness of a solution $u_0$ to \eqref{def-ellip} follows easily from the Lax-Milgram Lemma. Since $f+g \in L^2(0,T;V'_0)$, the estimate for the solution implies also $u_0\in L^2(0,T;V)$. 
\end{remark}

We shall sometimes use the corresponding problem with 
homogeneous boundary conditions:
\begin{equation}\label{weaksol2-ellip}
\begin{cases}
v_0(t)\in V_0 &\quad\text{for a.e.\ }t\in [0,T],\\
-\div(\A ev_0(t))=h(t)+\ell(t)&\quad\text{for a.e.\ }t\in [0,T],
\end{cases}
\end{equation}
with $h\in L^2(0,T;H)$ and $\ell\in H^1(0,T;V'_0)$.
\begin{remark}\label{z0}
The function $u_0$ is a solution to \eqref{def-ellip} if and only if $v_0=u_0-z$ is a solution to \eqref{weaksol2-ellip}
with
\begin{align*}
h(t)=f(t)\quad\text{and}\quad\ell(t)=g(t)+\div(\A ez(t)).
\end{align*}
\end{remark}

The following lemma will be used to prove the regularity with respect to time of the solution to~\eqref{def-ellip}.

\begin{lemma}\label{reg-ell-Hk}
Let $m\in \N$ and $p\in[1,+\infty)$. If $f=0$, $g\in W^{m,p}(0,T;V'_0)$, and $z\in W^{m,p}(0,T;V)$, then the solution $u_0$ to problem \eqref{def-ellip} satisfies $u_0\in W^{m,p}(0,T;V)$.
\end{lemma}

\begin{proof}
By Remark \ref{z0} it is enough to consider the case $z=0$. Let $R:V'_0\rightarrow V_0$ be the resolvent operator defined as follows:
\begin{equation*}
    R(\psi)=\varphi \iff \begin{cases}
\varphi\in V_0,\\
-\div(\A e\varphi)= \psi.
\end{cases}
\end{equation*}
Since $u_0(t)=R(g(t))$, the conclusion follows from the continuity of the linear operator $R$.
\end{proof}

\begin{remark}
In the case $f=0$, since $g\in W^{1,1}(0,T;V'_0)$ and $z\in W^{1,1}(0,T;V)$, we can apply Lemma \ref{reg-ell-Hk} to obtain that the solution $u_0$ to \eqref{def-ellip} belongs to $W^{1,1}(0,T;V)$, hence $u_0\in C^0([0,T];V)$.
\end{remark}

In the final statement of the next theorem, besides (H1)--(H3) we assume the following compatibility condition: there exists an extension of $g$ (still denoted by $g$) such that 
\begin{equation}\label{compat}
    g\in W^{1,1}(-a,T;V'_0)\quad \text{and}\quad -\div(\A e\uin(t))=g(t)\quad \text{for every $t\in [-a,0]$}.
\end{equation}

We are now in a position to state the main results of this paper. 

\begin{theorem}\label{main-thm-f+g}
Let us assume (H1)--(H3). Let $\u$ be the solution to the viscoelastic dynamic system \eqref{spaz} and let $u_0$ be the solution to the stationary problem \eqref{def-ellip}. Then 
\begin{alignat}{2}
\u&\xrightarrow[\varepsilon\to 0^+]{}u_0&&\qquad\text{strongly in $L^2(0,T;V)$}\label{conv-L2-f+g},\\
\varepsilon\du&\xrightarrow[\varepsilon\to 0^+]{}0&&\qquad\text{strongly in $L^2(0,T;H)$}\label{conv-L2d-f+g}.
\end{alignat}
If, in addition, $\fe=0$ for every $\varepsilon>0$, then
\begin{alignat}{2}
\u&\xrightarrow[\varepsilon\to 0^+]{}u_0&&\qquad\text{strongly in $L^{\infty}(\eta,T;V)$ for every $\eta\in(0,T)$},\label{th-4-g-f+g}\\
\varepsilon\du&\xrightarrow[\varepsilon\to 0^+]{}0&&\qquad\text{strongly in $L^{\infty}(\eta,T;H)$ for every $\eta\in(0,T)$}.\label{th-3-g-f+g}
\end{alignat}
If $\fe=0$ for every $\varepsilon>0$ and the compatibility condition \eqref{compat} holds, then we have also
\begin{alignat}{2}
\u&\xrightarrow[\varepsilon\to 0^+]{}u_0&&\qquad\text{strongly in $L^{\infty}(0,T;V)$}\label{th-4-g-inf},\\
    \varepsilon\du&\xrightarrow[\varepsilon\to 0^+]{}0&&\qquad\text{strongly in $L^{\infty}(0,T;H)$}.\label{conv-L2-der-inf}
\end{alignat}
\end{theorem}

In the case of solutions to problems \eqref{spaz2} we have the following results, assuming that
\begin{equation}\label{conv-dati-u}
    \uez\xrightarrow[\varepsilon\to 0^+]{}u^0\quad\text{strongly in $V$}\quad\text{and}\quad  \varepsilon\ueu\xrightarrow[\varepsilon\to 0^+]{}0\quad\text{strongly in $H$}.
\end{equation}

\begin{theorem}\label{main-f}
Let us assume (H1), (H2), and \eqref{conv-dati-u}. Let $\u$ be the solution to the viscoelastic dynamic system \eqref{spaz2}, with $\phie=\fe$ and $\gae=\gef$, and let $u_0$ be the solution to the stationary problem \eqref{def-ellip}. Then \eqref{conv-L2-f+g} and \eqref{conv-L2d-f+g} hold. Moreover, if $\fe=0$ for every $\varepsilon>0$, then \eqref{th-4-g-f+g} and \eqref{th-3-g-f+g} hold.
\end{theorem}

Theorems \ref{main-thm-f+g} and \ref{main-f} will be proved in several steps. First, we prove
\eqref{th-4-g-inf} and \eqref{conv-L2-der-inf} when $\fe=0$ and the compatibility condition \eqref{compat} holds (Theorem~\ref{teo-compat}). For $g\in H^2(0,T;V_0')$ the proof is based on the estimate in Lemma~\ref{lem:en-est} below, which is derived from the energy-dissipation balance \eqref{energy}. The general case is obtained by an approximation argument based on the same estimate.

Next, we prove that \eqref{conv-L2-f+g} holds for the solutions of \eqref{spaz2} if $\gae=\gamma=0$, $\ze=0$, $\uez=0$, and $\ueu=0$  (Proposition~\ref{new-main-th}).
The proof is obtained by means of a careful estimate of the solutions to the elliptic system \eqref{weak-fomulation-lap} obtained from \eqref{weaksol2} via Laplace Transform (Section \ref{laplace}). 
Under the general assumptions (H1), (H2), and \eqref{conv-dati-u} the same result is deduced from the previous one by an approximation argument based again on
Lemma~\ref{lem:en-est} below.

Then, \eqref{conv-L2d-f+g} is obtained from \eqref{conv-L2-f+g} using a suitable test function in \eqref{spaz2} (Theorem~\ref{uconvf+gder}). A further approximation argument gives \eqref{conv-L2-f+g} and \eqref{conv-L2d-f+g} under the assumptions (H1), (H2), and (H3) (Theorem~\ref{main-thm-f+g-part1}).

Finally, if $\fe=0$, we obtain \eqref{th-4-g-f+g} and \eqref{th-3-g-f+g} from \eqref{conv-L2-f+g} and \eqref{conv-L2d-f+g} (Lemma~\ref{regular-f+g}), concluding the proof of Theorems~\ref{main-thm-f+g} and \ref{main-f}.

The following lemma, derived from the energy-dissipation balance \eqref{energy}, will be frequently used to approximate the solutions of \eqref{weaksol2} by means of solutions corresponding to more regular data.

\begin{lemma}\label{lem:en-est}
Given $\varepsilon>0$, $\varphi_\varepsilon\in L^2(0,T;H)$, $\lee\in H^1(0,T;V'_0)$, $\vez\in V_0$, and $\veu\in H$, let $\v$ be the solution to \eqref{weaksol2} with $\he=\varepsilon\varphi_\varepsilon$. Then there exists a positive constant $C_{E}=C_{E}(\A,\B,\O,T)$, independent of $\varepsilon$, such that
\begin{align*}
    \varepsilon^2\|\dv\|^2_{L^{\infty}(0,T;H)}+\|\v\|^2_{L^{\infty}(0,T;V)}\leq C_{E}\Big(\varepsilon^2\|\veu\|^2+\|\vez\|^2_{V}+\|\varphi_\varepsilon\|^2_{L^1(0,T;H)}+\|\ell_\varepsilon\|^2_{W^{1,1}(0,T;V'_0)}\Big).
\end{align*}
\end{lemma}

\begin{proof}
By the energy-dissipation balance \eqref{energy} proved in Proposition \ref{diss-bal} and by \eqref{KP-ineq} and \eqref{CB2Q} there exists a positive constant $C=C(\A,\B,\O)$ such that 
\begin{equation}\label{en11}
    \varepsilon^2\|\dv(t)\|^2+\|\v(t)\|^2_V\leq C\Big(\varepsilon^2\|\veu\|^2+\|\vez\|_V^2+ \W_{\varepsilon}(t)\Big)\quad\text{for every $t\in [0,T]$},
\end{equation}
where the work is now defined by
\begin{equation}\label{en21}
  \W_{\varepsilon}(t)=\langle \lee(t),\v(t)\rangle-\langle \lee(0),\vez\rangle-\int_0^{t} \langle\dlee(\tau),\v(\tau)\rangle\de \tau+\int_0^t(\varphi_{\varepsilon}(\tau),\varepsilon\dv(\tau))\de\tau.
\end{equation}
Let $K_{\varepsilon}:=\varepsilon\|\dv(t)\|_{L^{\infty}(0,T;H)}$ and $E_{\varepsilon}:=\|\v(t)\|_{L^{\infty}(0,T;V)}$, which are finite by \eqref{cont}. Thanks to \eqref{en11} and \eqref{en21} for every $t\in [0,T]$ we get
\begin{align*}
    \varepsilon^2\|\dv(t)\|^2+\|\v(t)\|^2_V&\leq C\Big(\varepsilon^2\|\veu\|^2+\|\vez\|_V^2+\big(2\|\lee\|_{L^{\infty}(0,T;V'_0)}+\|\dlee\|_{L^1(0,T;V'_0)}\big)E_{\varepsilon}+\|\varphi_{\varepsilon}\|_{L^1(0,T;H)}K_{\varepsilon}\Big)\\
    &\leq C\Big(\varepsilon^2\|\veu\|^2+\|\vez\|_V^2+\big(3+\tfrac{2}{T}\big)\|\lee\|_{W^{1,1}(0,T;V'_0)}E_{\varepsilon}+\|\varphi_{\varepsilon}\|_{L^1(0,T;H)}K_{\varepsilon}\Big).
\end{align*}
By passing to the supremum with respect to $t$ and using the Young Inequality we can find a positive constant $C_{E}=C_{E}(\A,\B,\O,T)$ such that
\begin{equation*}
    K_{\varepsilon}^2+E_{\varepsilon}^2\leq C_{E}\Big(\varepsilon^2\|\veu\|^2+\|\vez\|_V^2+\|\varphi_{\varepsilon}\|^2_{L^1(0,T;H)}+\|\lee\|^2_{W^{1,1}(0,T;V'_0)}\Big),
\end{equation*}
which concludes the proof.
\end{proof}

In the proof of Theorem \ref{main-thm-f+g} we shall use the following lemma, which ensure that it is enough to consider the case $\ze=0$ and $z=0$.

\begin{lemma}\label{lemma-z}
If Theorem \ref{main-thm-f+g} holds when $\ze=0$ for every $\varepsilon>0$, then it holds for arbitrary $\{\ze\}_\varepsilon$ and $z$ satisfying (H2).
\end{lemma}

\begin{proof}
It is not restrictive to assume $\div(\B e\ze(0))=\div(\B ez(0))=0$. Indeed, if this is not the case, we can consider the solutions $\zez$ and $z^0$ to the stationary problems 
\begin{equation*}
\begin{cases}
\zez\in V_0,\\
-\div (\B e\zez)=\div(\B e\ze(0)),
\end{cases}
\quad\text{and}\quad
    \begin{cases}
z^0\in V_0,\\
-\div (\B ez^0)=\div(\B ez(0)),
\end{cases}
\end{equation*}  
and we can replace $\ze(t)$ and $z(t)$ by $\tze(t):=\ze(t)+\zez$ and $\tilde{z}(t):=z(t)+z^0$. It is clear that $\div(\B e\tilde{z}_{\varepsilon}(0))=\div(\B e\tilde z(0))=0$ and that problems \eqref{spaz} and \eqref{def-ellip} do not change passing from $\ze$ and $z$ to $\tze$ and~$\tilde{z}$. 

Let $\psi_{\varepsilon},\psi\colon[0,T]\rightarrow V'_0$ be the functions defined by
\begin{equation}\label{def-psi}
   \psi_{\varepsilon}(t):= \begin{cases}
0&\text{if } t\in(-\infty,0),\\
\div (\B e\ze(t))&\text{if } t\in[0,T],\\
\div (\B e\ze(T))&\text{if } t\in(T,+\infty),
\end{cases}
\quad\text{and}\quad
  \psi(t):= \begin{cases}
0&\text{if } t\in(-\infty,0),\\
\div (\B e z(t))&\text{if } t\in[0,T],\\
\div (\B e z(T))&\text{if } t\in(T,+\infty).
\end{cases}
\end{equation}
Since $\div(\B e\ze(0))=\div(\B e z(0))=0$, $\ze\in H^1(0,T;V)$, and $z\in  W^{1,1}(0,T;V)$, we obtain $\psi_{\varepsilon}\in H^1_{loc}(\R;V'_0)$ and $\psi\in W^{1,1}_{loc}(\R;V'_0)$. Moreover, thanks to (H2) we have 
\begin{equation}\label{convpsi}
    \psi_{\varepsilon}\xrightarrow[\varepsilon\to 0^+]{}\psi\quad\text{strongly in $W^{1,1}_{loc}(\R;V'_0)$}.
\end{equation}

Since $\u$ is the solution to \eqref{spaz}, by Remark \ref{corris} it solves \eqref{spaz2} with $\gae=\gef-p_{\varepsilon}$ and initial conditions defined by \eqref{u1}, where $p_{\varepsilon}$ is defined by \eqref{g0}. By Remark \ref{rem-um} the function $\v=\u-\ze$ is the solution to \eqref{weaksol2} with
\begin{equation}\label{def-gg}
\he(t)=\fe(t)-\varepsilon^2\ddze(t),\quad
    \lee(t)=\gef(t)-p_{\varepsilon}(t)+\div((\A+\B) e\ze(t))-\int_0^t\frac{1}{\beta\varepsilon}\e^{-\frac{t-\tau}{\beta\varepsilon}}\div(\B e\ze(\tau))\de\tau,
\end{equation}
and initial conditions $\vez$ and $\veu$ defined by \eqref{datin}. We define the family of convolution kernels $\{\rho_{\varepsilon}\}_{\varepsilon}\subset L^1(\R)$ by
\begin{equation}\label{rho}
    \rho_{\varepsilon}(t):=\begin{cases}
    \frac{1}{\beta\varepsilon}\e^{-\frac{t}{\beta\varepsilon}}&\text{if }t\in[0,+\infty),\\
    0&\text{if }t\in(-\infty,0),
    \end{cases} 
\end{equation}
and notice that, by \eqref{def-psi}, the integral in \eqref{def-gg} coincides with $(\rho_{\varepsilon}*\psi_{\varepsilon})(t)$, hence
\begin{equation*}
 \lee(t)=\gef(t)-p_{\varepsilon}(t)+\div(\A e\ze(t))+\psi_{\varepsilon}(t)-(\rho_{\varepsilon}*\psi_{\varepsilon})(t)\quad\text{for every $t\in[0,T]$.}
\end{equation*}

By Remark \ref{z0} the function $v_0=u_0-z$ is the solution to \eqref{weaksol2-ellip} with $h=f$ and $\ell=g+\div(\A ez)$. By the definition of $\v$ and $v_0$ it is clear that to prove the theorem it is enough to show that the conclusions of Theorem \ref{main-thm-f+g} holds for $\v$ and $v_0$. To this aim, we introduce the solution $\tve$ to \eqref{weaksol2} with $\he=\fe$, $\lee=\gef-p_{\varepsilon}+\div(\A e\ze)$, and $\vez$, $\veu$ defined by \eqref{datin}.  Then the function $\vb:=\v-\tve$ satisfies \eqref{weaksol2} with $\he=-\varepsilon^2\ddze$, $\lee=\psi_{\varepsilon}-\rho_{\varepsilon}*\psi_{\varepsilon}$, and homogeneous initial conditions. By Lemma \ref{lem:en-est} we can write
\begin{equation}\label{convy}
    \varepsilon^2\|\dvb\|^2_{L^{\infty}(0,T;H)}+\|\vb\|^2_{L^{\infty}(0,T;V)}\leq C_{E}\Big(\varepsilon^2\|\ddze\|^2_{L^1(0,T;H)}+\| \psi_{\varepsilon}-\rho_{\varepsilon}*\psi_{\varepsilon}\|^2_{W^{1,1}(0,T;V'_0)}\Big).
\end{equation}
By \eqref{convpsi} and by classical results on convolutions we obtain
\begin{equation*}
    \psi_{\varepsilon}-\rho_{\varepsilon}*\psi_{\varepsilon}\xrightarrow[\varepsilon\to 0^+]{} 0\quad\text{strongly in $W^{1,1}(0,T;V'_0)$}.
\end{equation*}
Since $\{\ddze\}_{\varepsilon}$ is bounded in $L^1(0,T;H)$ by (H2), from \eqref{convy} we deduce 
\begin{alignat}{2}
\v-\tve&\xrightarrow[\varepsilon\to 0^+]{}0&&\qquad\text{strongly in $L^{\infty}(0,T;V)$},\label{conv-med1*}\\
\varepsilon(\dv-\dot{\tilde v}_{\varepsilon})&\xrightarrow[\varepsilon\to 0^+]{}0&&\qquad\text{strongly in $L^{\infty}(0,T;H)$}.\label{conv-med2*}
\end{alignat}

By Remark \ref{corris} the function $\tve$ is the solution to \eqref{spaz} with $\gef$ replaced by $\gef+\div(\A e\ze)$ and $\ze=0$. Thanks to (H1) and (H2) we have
\begin{equation*}
    \gef+\div(\A e\ze)\xrightarrow[\varepsilon\to 0^+]{}g+\div(\A ez)\quad\text{strongly in $W^{1,1}(0,T;V'_0)$.}
\end{equation*} 
Since by hypothesis, Theorem \ref{main-thm-f+g} holds in the case of homogeneous boundary condition, its conclusions are valid for $\tve$ and $v_0$. Thanks to \eqref{conv-med1*} and \eqref{conv-med2*} the same results hold for $\v$ and $v_0$. This concludes the proof.
\end{proof}

In a similar way we can prove the following result.

\begin{lemma}\label{lemma-z2}
If Theorem \ref{main-f} holds when $\ze=0$ for every $\varepsilon>0$, then it holds for arbitrary $\{\ze\}_\varepsilon$ and $z$ satisfying (H2).
\end{lemma}

\section{The uniform convergence}

In this section we shall prove \eqref{th-4-g-inf} and \eqref{conv-L2-der-inf} of Theorem \ref{main-thm-f+g} under the compatibility condition \eqref{compat}. 

\begin{theorem}\label{teo-compat}
Let us assume (H1)--(H3), the compatibility condition \eqref{compat}, and $\fe=0$ for every $\varepsilon>0$. Let $\u$ be the solution to the viscoelastic dynamic system \eqref{spaz} and let $u_0$ be the solution to the stationary problem \eqref{def-ellip}, with $f=0$. Then \eqref{th-4-g-inf} and \eqref{conv-L2-der-inf} hold.
\end{theorem}

To prove the theorem we need the following lemma, which gives the result when $g$ is more regular.

\begin{lemma}\label{lem:appro}
Under the assumptions of Theorem \ref{teo-compat}, if $g\in H^2(0,T;V'_0)$, then \eqref{th-4-g-inf} and \eqref{conv-L2-der-inf} hold.
\end{lemma}
\begin{proof}
Thanks to Lemma \ref{lemma-z} we can suppose $z=0$ and $\ze=0$ for every $\varepsilon>0$. Let $p_{\varepsilon}$ be defined by \eqref{u1}. Since $\u$ is the solution to \eqref{spaz}, thanks to Remark \ref{corris} it solves \eqref{weaksol2} with $\he=0$, $\lee=\gef-p_{\varepsilon}$, $\vez=\uein(0)$, and $\veu=\duein(0)$. We fix $b>a>0$ and we extend the function $g$ in \eqref{compat} to $(-\infty,T)$ in such a way $g\in W^{1,1}(-\infty,T;V'_0)$ and $g(t)=0$ for every $t\in(-\infty,-b]$. Since $z=0$ we can extend $u_0$ by solving the following problem:
\begin{equation*}
\begin{cases}
u_0(t)\in V_0&\quad\text{for every $t\in (-\infty,T]$,}\\
-\div (\A eu_0(t))=g(t)&\quad\text{for every $t\in (-\infty,T]$}.
\end{cases}
\end{equation*}
We observe that $u_0=0$ on $(-\infty,-b]$ and $u_0=\uin$ on $[-a,0]$ by the compatibility condition \eqref{compat}. 

Assume $g\in H^2(0,T;V'_0)$. By Lemma \ref{reg-ell-Hk} (with $z=0$) we have $u_0\in H^2(0,T;V)$, hence by \eqref{def-ellip} we get
\begin{align}\label{iper-u0}
    \varepsilon^2\ddot u_0(t) &-\div((\A+\B) e u_0(t))+\int_0^t\frac{1}{\beta\varepsilon}\e^{-\frac{t-\tau}{\beta\varepsilon}}\div(\B e u_0(\tau))\de\tau\nonumber\\
    &= \varepsilon^2\ddot u_0(t)+g(t)-\div(\B e u_0(t))+(\rho_{\varepsilon}*\div(\B eu_0))(t)-\tilde{p}_{\varepsilon}(t)\quad\text{for a.e.\ }t\in [0,T],
\end{align}
where $\rho_{\varepsilon}$ is defined by \eqref{rho} and 
\begin{equation*}
    \tilde{p}_{\varepsilon}(t):=\e^{-\frac{t}{\beta\varepsilon}}\tilde{g}^0_{\varepsilon}\quad\text{with} \quad\tilde{g}^0_{\varepsilon}:=\int_{-\infty}^0\frac{1}{\beta\varepsilon}\e^{\frac{\tau}{\beta\varepsilon}}\div(\B eu_0(\tau))\de\tau=\int_{-b}^0\frac{1}{\beta\varepsilon}\e^{\frac{\tau}{\beta\varepsilon}}\div(\B eu_0(\tau))\de\tau.
\end{equation*}

Let  $q_{\varepsilon}:=\gef-g+\div(\B eu_0)-(\rho_{\varepsilon}*\div(\B eu_0))-p_{\varepsilon}+\tilde{p}_{\varepsilon}$.
By \eqref{iper-u0} the function $\ube:=\u-u_0$ satisfies \eqref{weaksol2} with $\he=-\varepsilon^2\ddot u_0$, $\lee=q_{\varepsilon}$, $\vez=\uein(0)-u_0(0)$, and $\veu=\duein(0)-\dot{u}_0(0)$.

Since $g\in W^{1,1}(-\infty,T;V'_0)$ and $g=0$ on $(-\infty,-b]$, by Lemma \ref{reg-ell-Hk} we obtain $u_0\in W^{1,1}(-\infty,T;V)$ and therefore $\div(\B eu_0)\in W^{1,1}(-\infty,T;V'_0)$. Then the properties of convolutions imply
\begin{equation}\label{rhoconv}
    \rho_{\varepsilon}*\div(\B eu_0)\xrightarrow[\varepsilon\to 0^+]{}\div(\B eu_0)\quad\text{strongly in $W^{1,1}(-\infty,T;V'_0)$}.
\end{equation}

As we have already observed, by the compatibility condition \eqref{compat} we have $u_0=\uin$ on $[-a,0]$, hence
\begin{align*}
    \|\tilde{g}^0_{\varepsilon}-g^0_{\varepsilon}\|_{V'_0}&\leq \int_{-\infty}^{-a}\frac{1}{\beta\varepsilon}\e^{\frac{\tau}{\beta\varepsilon}}\|\div(\B(e\uein(\tau))\|_{V'_0}\de\tau+\int_{-b}^{-a}\frac{1}{\beta\varepsilon}\e^{\frac{\tau}{\beta\varepsilon}}\|\div(\B(eu_0(\tau))\|_{V'_0}\de\tau\\
    &\hspace{2cm}+\|\div(\B (e\uein-e\uin))\|_{L^{\infty}(-a,0;V'_0)}.
\end{align*}
Thanks to (H3) we obtain $ \tilde{g}^0_{\varepsilon}-g^0_{\varepsilon}\rightarrow 0$ strongly in $V'_0$ as $\varepsilon\to 0^+$. Hence 
\begin{equation}\label{convlinf}
    \tilde{p}_{\varepsilon}-p_{\varepsilon}\xrightarrow[\varepsilon\to 0^+]{}0\quad\text{strongly in $W^{1,1}(0,T;V'_0)$.}
\end{equation}

By (H1), \eqref{rhoconv}, and \eqref{convlinf} we have
\begin{equation}\label{112}
    q_{\varepsilon}\xrightarrow[\varepsilon\to 0^+]{}0\quad\text{strongly in $W^{1,1}(0,T;V'_0)$}.
\end{equation}
Since $u_0(0)= \uin(0)$, (H3) gives
\begin{alignat}{2}
    \uein(0)-u_0(0)&\xrightarrow[\varepsilon\to 0^+]{}0&&\qquad\text{strongly in $V$},\label{113}\\
    \varepsilon(\duein(0)-\dot{u}_0(0))&\xrightarrow[\varepsilon\to 0^+]{}0&&\qquad\text{strongly in $H$}\label{114}.
\end{alignat}
By using Lemma \ref{lem:en-est} we get
\begin{align*}
     \varepsilon^2\|&\dube\|^2_{L^{\infty}(0,T;H)}+\| \ube\|_{L^{\infty}(0,T;V)}^2\nonumber\\
     &\leq C_{E}\Big(\varepsilon^2\|\duein(0)-\dot{u}_0(0)\|^2+\| \uein(0)-u_0(0)\|^2_{V}+\varepsilon^2\|\ddot u_0\|^2_{L^1(0,T;H)}+\|q_{\varepsilon}\|^2_{W^{1,1}(0,T;V'_0)}\Big),
\end{align*}
therefore thanks to \eqref{112}, \eqref{113}, and \eqref{114} we obtain the conclusion.
\end{proof}

In the proof of Theorems \ref{teo-compat}, \ref{uconvf+g}, and \ref{main-thm-f+g-part1} we shall use the following density result.

\begin{lemma}\label{densita}
Let $X,Y$ be two Hilbert spaces such that $X\hookrightarrow Y$ continuously, with $X$ dense in $Y$. Then for every $m,n\in\N$ with $m\leq n$, and $p\in [1,2]$ the space $H^n(0,T;X)$ is dense in $W^{m,p}(0,T;Y)$.
\end{lemma}

\begin{proof}
Since every simple function with values in $Y$ can be approximated by simple functions with values in $X$, it is easy to see that $L^2(0,T;X)$ is dense in $L^p(0,T;Y)$.

To prove the result for $m=1$ we fix $u\in W^{1,p}(0,T;Y)$. By the density of $L^2(0,T;X)$ in $L^p(0,T;Y)$ we can find a sequence $\{\psi_k\}_k\subset L^2(0,T;X)$ such that $\psi_k\rightarrow\dot{u}$ strongly in $L^p(0,T;Y)$ as $k\to+\infty$.
By the density of $X$ in $Y$ there exists $\{u^0_k\}_k\subset X$ such that $u_k^0\rightarrow u(0)$ strongly in $Y$ as $k\to +\infty$. Now we define 
\begin{equation*}
    u_k(t):=\int_0^t\psi_k(\tau)\de \tau+u_k^0.
\end{equation*}
It is easy to see that $\{u_k\}_k\subset H^1(0,T;X)$ and $u_k\rightarrow u$ strongly in $W^{1,p}(0,T;Y)$ as $k\to +\infty$.

Arguing by induction we can prove that for every integer $m\geq 0$ the space $H^m(0,T;X)$ is dense in $W^{m,p}(0,T;Y)$. Since $H^n(0,T;X)$ is dense in $H^m(0,T;X)$, the conclusion follows.
\end{proof}

We are now in a position to deduce Theorem \ref{teo-compat} from Lemma \ref{lem:appro} by means of an approximation argument.

\begin{proof}[Proof of Theorem \ref{teo-compat}] Thanks to Lemma \ref{lemma-z} we can suppose $z=0$ and $\ze=0$ for every $\varepsilon>0$. We fix $\delta>0$. By Lemma \ref{densita} there exists a function $\psi\in H^2(0,T;V'_0)$ such that 
\begin{equation}\label{eta1}
\|\psi-g\|_{W^{1,1}(0,T;V'_0)}<\delta.    
\end{equation}
By (H1) there exists a positive number $\varepsilon_0=\varepsilon_0(\delta)$ such that
\begin{equation}\label{eta2}
\|\psi-\gef\|_{W^{1,1}(0,T;V'_0)}<\delta\qquad\text{for every $\varepsilon\in(0,\varepsilon_0)$}.    
\end{equation}

Let $p_{\varepsilon}$ be defined by \eqref{g0}. Since $\u$ is the solution to \eqref{spaz} with $\fe=0$ and $\ze=0$, thanks to Remark \ref{corris} it solves \eqref{weaksol2} with $\he=0$, $\lee=\gef-p_{\varepsilon}$, $\vez=\uein(0)$, and $\veu=\duein(0)$. Moreover, let $\tue$ be solution to \eqref{weaksol2} with $\he=0$, $\lee=\psi-p_{\varepsilon}$, $\vez=\uein(0)$, and $\veu=\duein(0)$, and let $\tilde{u}_0$ be the solution to \eqref{weaksol2-ellip} with $h=0$ and $\ell=\psi$. Thanks to Remark \ref{corris} the function $\tue$ is the solution to \eqref{spaz} with $\fe=0$, $\gef=\psi$, and $\ze=0$, hence by Lemma \ref{lem:appro} we have
\begin{alignat}{2}
    \tue&\xrightarrow[\varepsilon\to 0^+]{}\tilde{u}_0&&\qquad\text{strongly in $L^{\infty}(0,T;V)$},\label{conv1t}\\
   \varepsilon\dtue&\xrightarrow[\varepsilon\to 0^+]{}0&&\qquad\text{strongly in $L^{\infty}(0,T;H)$}.\label{conv2t}
\end{alignat}

We now consider the functions $\bar{u}_0:=\tilde{u}_0-u_0$ and $\ube:=\tue-\u$. Since $\bar{u}_0$ is the solution to \eqref{weaksol2-ellip}, with $h=0$ and $\ell=\psi-g$, by the Lax-Milgram Lemma we get
\begin{equation}\label{stima-u0}
    \|\bar{u}_0\|_{L^{\infty}(0,T;V)}\leq \tfrac{C^2_P+1}{c_{\A}} \|\psi-g\|_{L^{\infty}(0,T;V'_0)}\leq \tfrac{C^2_P+1}{c_{\A}}(1+\tfrac{1}{T}) \|\psi-g\|_{W^{1,1}(0,T;V'_0)}.
\end{equation}
Moreover, since $\ube$ is the solution to \eqref{weaksol2}, with $\he=0$, $\lee=\psi-\gef$, $\vez=0$, and $\veu=0$, thanks to Lemma \ref{lem:en-est} we get
\begin{equation}\label{en-es-2}
    \varepsilon^2\|\dube\|^2_{L^{\infty}(0,T;H)}+\|\ube\|^2_{L^{\infty}(0,T;V)}\leq C_E\|\psi-\gef\|^2_{W^{1,1}(0,T;V'_0)}.
\end{equation}
By combining \eqref{eta1}, \eqref{eta2}, \eqref{stima-u0}, and \eqref{en-es-2}, we can find a positive constant $C=C(\A,\B,\O,T)$ such that
\begin{equation}\label{en-es-2*}
    \varepsilon\|\dube\|_{L^{\infty}(0,T;H)}+\|\ube\|_{L^{\infty}(0,T;V)}+\|\bar{u}_0\|_{L^{\infty}(0,T;V)}\leq  C\delta\quad\text{for every $\varepsilon\in(0,\varepsilon_0)$}.
\end{equation}
Since 
\begin{align*}
    \|\u-u_0\|_{L^{\infty}(0,T;V)}&\leq \|\ube\|_{L^{\infty}(0,T;V)}+\|\tue-\tilde{u}_0\|_{L^{\infty}(0,T;V)}+\|\bar{u}_0\|_{L^{\infty}(0,T;V)},\\
    \varepsilon\|\du\|_{L^{\infty}(0,T;H)}&\leq \varepsilon\|\dube\|_{L^{\infty}(0,T;H)}+\varepsilon\|\dtue\|_{L^{\infty}(0,T;H)},
\end{align*}
by \eqref{conv1t}, \eqref{conv2t}, and \eqref{en-es-2*} we have
\begin{align*}
    \limsup_{\varepsilon\to 0^+}\|\u-u_0\|_{L^{\infty}(0,T;V)}\leq C\delta\quad\text{and}\quad\limsup_{\varepsilon\to 0^+}\|\varepsilon\du\|_{L^{\infty}(0,T;H)}\leq C\delta.
\end{align*}
The conclusion follows from the arbitrariness of $\delta>0$.
\end{proof}

\section{Use of the Laplace Transform}\label{laplace}

In this section we shall use the Laplace Transform to prepare the proof of 
the convergence, as $\varepsilon\to 0^+$, of the solutions of the problems
\begin{equation}\label{weak-form-v}
\begin{cases}
\v\in \VD\\
\displaystyle \varepsilon^2\ddv(t) -\div((\A+\B) e\v(t))+\int_0^t\frac{1}{\beta\varepsilon}\e^{-\frac{t-\tau}{\beta\varepsilon}}\div(\B e\v(\tau)) \de\tau= \he(t)\quad\text{for a.e.\ }t\in [0,T],\\
\v(0)=0\quad\text{in $H$}\quad \text{and}\quad \dv(0)=0\quad\text{in }V'_0,
\end{cases}
\end{equation}
to the solution $v_0$ of the problem
\begin{equation}
\begin{cases}
v_0(t)\in V_0&\quad\text{for a.e.\ }t\in [0,T],\\
 -\div(\A ev_0(t))=h(t)&\quad\text{for a.e.\ }t\in [0,T],\label{quas}
\end{cases}
\end{equation}
when $\{\he\}_{\varepsilon}\subset L^2(0,T;H)$, $h\in L^2(0,T;H)$, and
\begin{equation}\label{conv-fe}
    \he\xrightarrow[\varepsilon\to 0^+]{}h\qquad \text{strongly in }L^2(0,T;H),
\end{equation}
This partial result will be the starting point for the proof of the convergence in $L^2(0,T;V)$ under the general assumptions of Theorem \ref{main-thm-f+g}. 

\subsection{The Laplace Transform for functions with values in Hilbert spaces}
Given a complex Hilbert space $X$, let $r\in L^1_{loc}(0,+\infty;X)$ be a function such that
\begin{equation}\label{LT}
   \int_0^{+\infty}\e^{-\alpha t}\|r(t)\|_X\,\de t<+\infty\quad \text{for every $\alpha>0$,}
\end{equation}
and let $\C_+:=\{s\in\C:\Re(s)>0\}$. The Laplace Transform of $r$ is the function $\hat r:\C_+\rightarrow X$ defined by
\begin{equation}\label{Lap-Tran}
    \hat{r}(s):=\int_0^{+\infty}\e^{-st}r(t)\de t\quad \text{for every $s\in\C_+$}.
\end{equation}
Besides $\hat r$, we shall also use the notation $\mathcal{L}(r)$, which is sometimes written as $\mathcal{L}_t(r(t))$, with dummy variable~$t$. In the particular case $r\in L^{\infty}(0,+\infty;X)$ we have
\begin{equation*}
    \|\hat{r}(s)\|_X\leq \frac{1}{s_1}\|r\|_{L^{\infty}(0,+\infty;X)}\quad \text{for every $s=s_1+is_2\in \C_+$, with $s_1,s_2\in\R$}.
\end{equation*}

There is a close connection between the Laplace Transform and the Fourier Transform, defined for every $\rho\in L^{1}(\R;X)$ as the function $\mathcal{F}(\rho)\in L^{\infty}(\R;X)$ given by
\begin{equation}\label{Fou-Tran}
    \mathcal{F}(\rho)(\xi)=\int_{-\infty}^{+\infty}\e^{-i\xi t}\rho(t)\de t\quad \text{for every $\xi\in\R$}.
\end{equation}
For $\mathcal{F}(\rho)$ we use also the notation $\mathcal{F}_t(\rho(t))$ with dummy variable $t$. For the main properties of the Fourier and Laplace Transforms of functions with values in Hilbert spaces we refer to \cite{LT}.

We extend the function $r$ satisfying \eqref{LT} by setting $r(t)=0$ for every $t<0$. By \eqref{Lap-Tran} and \eqref{Fou-Tran} we have
\begin{equation*}
    \mathcal{L}_t(r(t))(s)=\mathcal{F}_{t}(\e^{-s_1t} r(t))(s_2)\quad \text{for every $s=s_1+is_2\in \C_+$, with $s_1,s_2\in\R$}.
\end{equation*}

We remark that the Laplace Transform commutes with linear transformations, as shown in the following proposition (see, for instance \cite[Proposition 1.6.2]{LT}).
\begin{proposition}\label{prop:lin-laplace}
Let $X$ and $Y$ be two complex Hilbert spaces, let $r \in L^1_{loc}(0,+\infty;X)$, and let $T$ be a continuous linear operator from $X$ to $Y$. Then $T\circ r\in L^1_{loc}(0,+\infty;Y)$. If in addition, $r$ satisfies \eqref{LT}, then the same property holds also for $T\circ r$, with $X$ replaced by $Y$, and $\mathcal{L}(T\circ r)(s)=(T\circ\hat{r})(s)$ for every $s\in \C_+$.
\end{proposition}

Now we consider the Inverse Laplace Transform. Let $R:\C_+\rightarrow X$ be a function. Suppose that there exists $r\in L^1_{loc}(0,+\infty;X)$ such that \eqref{LT} holds and $\mathcal{L}(r)(s)=R(s)$ for every $s\in\C_+$. In this case we say that $r$ is the Inverse Laplace Transform of $R$, and we use the notation $r=\mathcal{L}^{-1}(R)$ or $r=\mathcal{L}_s^{-1}(R(s))$ with dummy variable $s$. It can be proven that $r$ is uniquely determined up to a negligible set (see \cite[Theorem 1.7.3]{LT}). Moreover, $r$ can be obtained by the Bromwich Integral Formula:
\begin{equation}\label{lap-inv}
    r(t)=\mathcal{L}^{-1}(R)(t)=\frac{\e^{s_1 t}}{2\pi }\lim_{k\to +\infty }\int_{-k}^{k}e^{is_2 t}R(s_1+is_2) \de s_2,
\end{equation}
where $s_1$ is an arbitrary positive number. Clearly \eqref{lap-inv} can be expressed in terms of the Inverse Fourier Transform, namely
\begin{equation}\label{inv-lap-four}
    r(t)=\mathcal{L}_s^{-1}(R(s))(t)=\e^{s_1 t}\mathcal{F}_{s_2}^{-1}(R(s_1+is_2))(t),
\end{equation}
where $\mathcal{F}^{-1}_{s_2}(R(s_1+is_2))$ denotes the Inverse Fourier Transform with respect to the variable $s_2$.

To use the Laplace Transform, we extend our problems from the interval $[0,T]$ to $[0,+\infty)$. To do this, we extend the functions $\he$ and $h$, introduced in \eqref{conv-fe}, by setting them equal to zero in $(T,+\infty)$, and we consider the solution to \eqref{weak-form-v} in $[0,+\infty)$, which we still denote $\v$. Moreover, we consider the solution to \eqref{quas} in $[0,+\infty)$, which we still denote $v_0$. Notice that, thanks to the choice of the extension we have
\begin{equation*}
    \he\xrightarrow[\varepsilon\to 0^+]{}h\qquad \text{strongly $L^2(0,+\infty;H)$}.
\end{equation*}

By Proposition \ref{lem:en-est} and by using the equality $\he=0$ on $(T,+\infty)$, we get
\begin{equation}\label{cont-v}
    \v\in  L^{\infty}(0,+\infty;V_0)\quad\text{and}\quad  \dv\in  L^{\infty}(0,+\infty;H).
\end{equation}
Since $h\in L^2(0,T;H)$ and $h=0$ on $(T,+\infty)$, by means of standard estimates for the solution to \eqref{quas} we obtain 
\begin{equation}\label{v=0}
    v_0\in L^2(0,+\infty;V_0)\quad\text{and}\quad  v_0=0 \quad\text{ on $(T,+\infty)$}.
\end{equation}
From \eqref{CB1Q}, \eqref{weak-form-v}, and \eqref{cont-v} we can deduce
\begin{equation}\label{v-sec}
\ddv\in L^2(0,T;V_0')\cap L^{\infty}(T,+\infty;V_0').    
\end{equation}

To study our problem by means of the Laplace Transform we introduce the complexification of the Hilbert spaces $H$, $V_0$, and $V'_0$ defined by
\begin{equation*}
    \hat{H}:=L^2(\O;\C^d),\quad \hat{V}_0:=H^1(\O;\C^d),\quad \hat{V}_0':=H^{-1}(\O;\C^d).
\end{equation*}
The symbols $(\cdot,\cdot)$ and $\|\cdot\|$ denote the hermitian product and the norm in $\hat{H}$ or in other complex $L^2$ spaces. For every $s\in\C_+$ the Laplace Transforms $\hat{h}_{\varepsilon}(s)$ and $\hat{h}(s)$ of $h_{\varepsilon}$ and $h$ in $\hat{H}$ are well defined. Thanks to \eqref{cont-v} and \eqref{v=0} the Laplace Transforms $\hv(s)$ and $\hat{v}_0(s)$ in $\hat{V}_0$ are well defined for every $s\in\C_+$. By \eqref{v-sec} the Laplace Transform $ \hat{\ddot{v}}_{\varepsilon}(s)$ of $\ddv(s)$ in $\hat{V}_0'$ is well defined for every $s\in\C_+$. Using \eqref{cont-v} we can integrate by parts two times in the integral which defines $ \hat{\ddot{v}}_{\varepsilon}$ and, since $\v(0)=0$ and $\dv(0)=0$, we  obtain 
\begin{equation}\label{rel-v-ddv}
   \hat{\ddot{v}}_{\varepsilon}(s)=s^2\hv(s)\quad \text{for every $s\in \C_+$}.
\end{equation}

By considering the operators $S_{\A},S_{\B}:\hat{V}_0\rightarrow \hat{V}_0'$ defined by
\begin{align*}
    S_{\A}(\psi):=-\div(\A e\psi)\quad\text{and}\quad  S_{\B}(\psi):=-\div(\B e\psi),
\end{align*}
we can rephrase \eqref{weak-form-v} and \eqref{quas} as equalities of elements of $\hat{V}_0'$: 
\begin{equation}\label{new-form-v}
   \varepsilon^2 \ddv(t) = S_{\B}\Big(\int_0^t\frac{1}{\beta \varepsilon}\e^{-\frac{t-\tau}{\beta\varepsilon}} \v(\tau) \de \tau\Big)- (S_{\A}+S_{\B})(\v(t))+\he(t)\quad\text{for a.e.\  $t\in [0,+\infty)$},
\end{equation}
\begin{equation}\label{new-form-ell}
    S_{\A}(v_0(t))=h(t)\quad\text{for a.e.\  $t\in [0,+\infty)$}.
\end{equation}

Now we want to consider the Laplace Transforms, in the sense of $\hat{V}'_{0}$, of both sides of these equations. By Proposition \ref{prop:lin-laplace} we can say 
\begin{equation}\label{S_2}
\mathcal{L}(S_{\A}(\v))=S_{\A}(\hv), \quad   \mathcal{L}(S_{\B}(\v))=S_{\B}(\hv),\quad \mathcal{L}(S_{\A}(v_0))=S_{\A}(\hat{v}_0),
\end{equation}
where $\hv$ and $\hat{v}_0$ are the Laplace Transforms of $\v$ and $v_0$, respectively, in the sense of $\hat{V}_0$.
Moreover, since we have
\begin{equation*}
     \sup_{t\in[0,+\infty)}\Big\|\int_0^t\frac{1}{\beta\varepsilon}\e^{-\frac{t-\tau}{\beta\varepsilon}} \v (\tau)\de \tau\Big\|_{V_0}\leq  \| \v\|_{L^{\infty}(0,+\infty;V_0)},
\end{equation*}
this integral admits Laplace Transform in the sense of $\hat{V}_0$, which for every $s\in\C_+$ satisfies 
\begin{align*}
    \mathcal{L}_t\Big(\int_0^{t}\frac{1}{\beta\varepsilon}\e^{-\frac{t-\tau}{\beta\varepsilon}}\v (\tau)\de \tau\Big)(s)=\frac{1}{\beta\varepsilon s+1}\hv(s).
\end{align*}
Hence, by using Proposition \ref{prop:lin-laplace} again, we obtain
\begin{equation}\label{Lap-Tra-2}
\mathcal{L}_t\Big(S_{\B}\Big(\int_0^t\frac{1}{\beta\varepsilon}\e^{-\frac{t-\tau}{\beta\varepsilon}} \v (\tau)\de \tau\Big)\Big)(s)=\frac{1}{\beta\varepsilon s+1}S_{\B}(\hv(s)).
\end{equation}

\subsection{Properties of the Laplace Transform of the solutions}
Thanks to \eqref{rel-v-ddv}, \eqref{S_2}, and \eqref{Lap-Tra-2} we can apply the Laplace Transform to both sides of \eqref{new-form-v} and \eqref{new-form-ell} to deduce the following equalities in $\hat{V}_0'$:
\begin{equation}\label{weak-fomulation-lap}
  \varepsilon^2 s^2\hv(s)-\div((\A+\B)e\hv(s))+\frac{1}{\beta\varepsilon s+1}\div(\B e\hv(s))=\hh(s)\quad\text{for every $s\in\C_+$},
\end{equation}
\begin{equation}\label{weak-fomulation-ellip-lap}
    -\div(\A e\hat{v}_0(s))=\hat{h}(s)\quad\text{for every $s\in\C_+$}.
\end{equation}

Our purpose is to prove that for every $s_1>0$ we have 
\begin{equation}\label{dadim}
    \int_{-\infty}^{+\infty}\|\hv(s_1+is_2)-\hat{v}_0(s_1+is_2)\|_{\hat{V}_0}^2\de s_2\xrightarrow[\varepsilon\to 0^+]{} 0.
\end{equation} 

To prove \eqref{dadim} we need two lemmas. In the first one we deduce from \eqref{weak-fomulation-lap} an estimate for $\hv(s)$, which is used in the second lemma to prove a convergence result for $\hv(s)$. 

\begin{lemma}\label{lem:ex-lax}
For every $s\in \C_+$ there exists a positive constant $M(s)$ such that
\begin{equation}\label{lax-estimate}
    \|\hv(s)\|_{\hat{V}_0}\leq M(s)\|\hh(s)\|\quad \text{for every $\varepsilon\in(0,1)$}.
\end{equation}
\end{lemma}

\begin{proof}
We fix $\varepsilon\in(0,1)$ and for every $s\in \C_+$ we define the operator $S_{\varepsilon}(s)\colon \hat{V}_0\rightarrow \hat{V}_0'$ in the following way:
\begin{equation*}
     S_{\varepsilon}(s)(\psi):=\varepsilon^2 s^2\psi-\div((\A+\B)e\psi)+\frac{1}{\beta\varepsilon s+1}\div(\B e\psi)\qquad \text{for every $\psi\in \hat{V}_0$}.
\end{equation*}
Since $S_{\varepsilon}(s)(\hv(s))=\hh(s)$ by \eqref{weak-fomulation-lap}, the Lax-Milgram Lemma, together with the Korn-Poincaré Inequality \eqref{KP-ineq}, implies \eqref{lax-estimate} if we can show that for every $s\in \C_+$ there exists a positive constant $K(s)$, independent of $\varepsilon$, such that 
\begin{align}\label{lax-2}
    c_{\A}K(s)\|e\psi\|^2\leq |\langle S_{\varepsilon}(s)(\psi),\psi\rangle| =\frac{|(\beta\varepsilon^3 s^3+\varepsilon^2 s^2)\|\psi\|^2+\beta\varepsilon s((\A+\B)e\psi,e\psi)+(\A e\psi,e\psi)|}{|\beta\varepsilon s +1|}
\end{align}
for every $\psi\in \hat{V}_0$. 

We can suppose $\psi\in \hat{V}_0\setminus\{0\}$, otherwise the inequality is trivially satisfied, and we set 
\begin{equation*}
    a:=\frac{(\A e\psi,e\psi)}{\|\psi\|^2} \quad \text{and}\quad b:=\frac{((\A+\B)e\psi,e\psi)}{\|\psi\|^2},
\end{equation*}
which satisfy, thanks to the Korn-Poincaré Inequality \eqref{KP-ineq} and to \eqref{CB1Q}--\eqref{CB2Q}, the following relations 
\begin{align}
    &a\geq \frac{c_{\A}\|e\psi\|^2}{\|\psi\|^2}\geq \frac{c_{\A}}{C_P^2}=:a_0, \qquad b\geq \frac{(c_{\A}+c_{\B})\|e\psi\|^2}{\|\psi\|^2}\geq \frac{c_{\A}+c_{\B}}{C_P^2}=:b_0,\qquad a \leq c_0a\leq b\leq c_1a,\label{coercQ}
\end{align}
where $c_0:=1+\frac{c_{\B}}{C_{\A}}$ and $c_1:=1+\frac{C_{\B}}{c_{\A}}$.
Therefore, to prove \eqref{lax-2} it is enough to obtain
\begin{align}\label{lax}
    \left|\frac{\beta \varepsilon^3 s^3+\varepsilon^2 s^2+\beta b \varepsilon s +a}{\beta\varepsilon s +1}\right|\geq K(s)\, a\quad \text{for every $s\in\C_+$}.
\end{align}

For simplicity of notation we set $z=\varepsilon s$ and we consider two cases.\\
\noindent {\it Case $b>\frac{2}{3\beta^2}$}. In this situation, thanks to \eqref{sistema} we know that the polynomial $\beta z^3+z^2+\beta b z+a$ has one real root $z_0$ and two complex and conjugate ones $w$ and $\Bar{w}$. Therefore, thanks to Lemmas \ref{comp-bound}  and \ref{comp-in}, we can write
\begin{align}\label{dis-co}
    \left|\frac{\beta z^3+z^2+\beta b z+a}{\beta z+1}\right|&=\left|\frac{\beta (z-z_0)(z-w)(z-\Bar{w})}{\beta z+1}\right|\geq \left|\frac{\beta(z-z_0)}{\beta z+1}\right||\Re(w)||\Im(w)|\nonumber\\
    &=\left|\frac{\beta(z-z_0)}{\beta z+1}\right||\Re(w)|\sqrt{3|\Re(w)|^2+\frac2\beta\Re(w)+b}\geq \left|\frac{\beta(z-z_0)}{\beta z+1}\right|\alpha\sqrt{b-\frac{1}{3\beta^2}}\nonumber\\
    &\geq \left|\frac{\beta(z-z_0)}{\beta z+1}\right|\alpha\sqrt{\frac{ b}{2}}\geq 
    \frac\alpha{\sqrt{3}}\left|\frac{z}{\beta z+1}\right|,
\end{align}
where in the last inequality we used $z_0<0$.

If $a\leq 2|z|^2$, then $|z|\geq \frac{a}{2|z|}$ and, thanks to \eqref{dis-co}, we deduce
\begin{equation}\label{DM-A}
    \left|\frac{\beta z^3+z^2+\beta b z+a}{\beta z+1}\right|\geq \frac{\alpha}{2\sqrt{3}}\frac{a}{|z(\beta z+1)|}.
\end{equation}
For $a>2|z|^2$ we have
\begin{equation*}
     \frac{1}{a}\left|\frac{\beta z^3+z^2+\beta b z+a}{\beta z+1}\right|=\left|\frac{ z^2}{a}+\frac{\beta bz+a}{a(\beta z+1)}\right|\geq \left|\frac{\beta bz+a}{a(\beta z+1)} \right|-\frac{1}{2},
\end{equation*}
and, by writing $z=x+iy$, we obtain
\begin{equation*}
    \left|\frac{\beta bz+a}{a(\beta z+1)}   \right|=\left|\frac{\beta b x+a+i\beta b y}{\beta a x+a+i\beta ay}   \right|=\sqrt{\frac{(\beta b x+a)^2+\beta^2 b^2y^2}{(\beta a x+a)^2+\beta^2a^2y^2}}\geq 1,
\end{equation*}
 which implies
\begin{equation}\label{DM-B}
   \left|\frac{\beta z^3+z^2+\beta bz+a}{\beta z+1}\right|\geq \frac{a}{2}.
\end{equation}
By \eqref{DM-A} and \eqref{DM-B} in the case $b> \frac{2}{3\beta^2}$ we conclude
\begin{equation}\label{DM-C}
    \left|\frac{\beta z^3+z^2+\beta b z+a}{\beta z+1}\right|\geq
    \min\Big\{\frac12,\frac{\alpha}{2\sqrt{3}}\frac{1}{|z(\beta z+1)|}\Big\}a.
\end{equation}
\noindent {\it Case $b_0\leq b\leq \frac{2}{3\beta^2}$.} In this case, thanks to \eqref{coercQ}, we have $a_0\leq a\leq \frac{2}{3\beta^2}$. We define
\begin{equation*}
    R:=\sqrt{\frac{2(2+c_1)}{3\beta^2}}.
\end{equation*}
Then for $z\in\C_+$, with $|z|>R$, we get
\begin{equation}\label{DM*}
    \frac{1}{a}\left|\frac{\beta z^3+z^2+\beta b z+a}{\beta z+1}\right|=\left|\frac{ z^2}{a}+\frac{\beta bz+a}{a(\beta z+1)}\right|\geq \frac{3\beta^2|z|^2}{2}-\frac{b}{a}\left|\frac{\beta z}{\beta z+1}\right|-\frac{1}{|\beta z+1|}\geq 2+c_1-c_1-1=1,
\end{equation}
where we used the inequalities $|\beta z|\leq |\beta z+1|$ and $1\leq |\beta z+1|$.

To deal with the case $z\in \C_+$, with $|z|\leq R$, we define
\begin{equation*}
   \gamma:=\min\left\{\left|\frac{\beta z^3+z^2+\beta b z+a}{a(\beta z+1)}\right|:\Re(z)\geq 0,\quad |z|\leq R,\quad b_0\leq b\leq \frac{2}{3\beta^2}, \quad a_0\leq a\leq \frac{2}{3\beta^2 } \right\},
\end{equation*}
and we claim $\gamma>0$. Indeed the function under examination is continuous with respect to $(z, a, b)$, and
by Lemma \ref{comp-bound} it does not vanish in the compact set considered in the minimum problem.
By using also \eqref{DM*} we conclude that for $b_0\leq b\leq \frac{2}{3\beta^2}$ we have
\begin{equation}\label{b-little}
     \left|\frac{\beta z^3+z^2+\beta b z +a}{\beta z+1}\right|\geq \min\{\gamma,1\} a.
\end{equation}
for every $z\in \C_+$ and every $a$ satisfying \eqref{coercQ}.
Since $\varepsilon\in (0,1)$ we have
\begin{equation*}
    \frac{1}{|\varepsilon s(\beta \varepsilon s+1)|}\geq\frac{1}{| s(\beta  s+1)|},
\end{equation*}
therefore, by setting
\begin{equation*}
    K(s):=
    \min\Big\{\frac12,\frac{\alpha}{2\sqrt{3}}\frac{1}{|s(\beta s+1)|},\gamma\Big\},
\end{equation*}
from \eqref{DM-C} and \eqref{b-little} we obtain \eqref{lax}, which concludes the proof. 
\end{proof}

\subsection{Convergence of the Laplace Transform of the solutions}
We begin by proving the pointwise convergence.

\begin{lemma}\label{lem:conv-punt-hv}
For every $s\in\C_+$ we have
\begin{equation*}
    \hv(s)\xrightarrow[\varepsilon\to 0^+]{}\hat{v}_0(s)\qquad \text{strongly in $\hat{V}_0$}.
\end{equation*}
\end{lemma}

\begin{proof}
Thanks to \eqref{conv-fe} and to the Hölder Inequality for every $s\in \C_+$ we get 
\begin{align}\label{punt-G-conv}
     \|\hh(s)-\hat{h}(s)\|\leq \int_0^{+\infty}\e^{-\Re(s)t}\|\he(t)-h(t)\|\de t\leq \frac{1}{\sqrt{2\Re(s)}}\|\he-h\|_{L^2(0,T;H)}\xrightarrow[\varepsilon\to 0^+]{}0.
\end{align}
Consequently, thanks to Lemma \ref{lem:ex-lax}, for every $s\in\C_+$ there exist two constants $\Bar{M}(s)>0$ and $\varepsilon(s)\in(0,1)$ such that
\begin{equation}\label{bound-hv}
    \|\hv(s)\|_{\hat{V}_0}\leq \Bar{M}(s)\qquad \text{for every $\varepsilon\in (0,\varepsilon(s))$}.
\end{equation}
By \eqref{bound-hv} we can say that for every $s\in\C_+$ there exist a sequence $\varepsilon_j\xrightarrow[]{}0^+$ and $v^*(s)\in \hat{V}_0$ such that  
\begin{equation}\label{weak-conv-hv}
    \hvj(s)\xrightharpoonup[j\to +\infty]{}v^*(s)\qquad\text{weakly in $\hat{V}_0$}.
\end{equation}
Thanks to \eqref{sym} and \eqref{weak-conv-hv} for every $\psi\in \hat{V}_0$ we deduce
\begin{align*}
    &((\A+\B)e\hvj(s),e\psi)\xrightarrow[j\to +\infty]{}((\A+\B)ev^*(s),e\psi),\quad |\varepsilon_j^2s^2(\hvj(s),\psi)|\leq \varepsilon_j^2|s|^2\Bar{M}(s)\|\psi\|\xrightarrow[j\to +\infty]{}0,\\
    &\begin{aligned}
        \Big|\frac{1}{\beta \varepsilon_j s+1}(\B e\hvj(s),e\psi)-(\B ev^*(s),e\psi)\Big|&\leq\left|(\B (e\hvj(s)-ev^*(s)),e\psi)\right|+\frac{\beta\varepsilon_j|s|}{|\beta \varepsilon_j s+1|}|(\B e\hvj(s),e\psi)|\\
    &\leq \left|( e\hvj(s)-ev^*(s),\B e\psi)\right|+\beta\varepsilon_j|s|C_{\B}\Bar{M}(s)\|e\psi\|\xrightarrow[j\to +\infty]{}0.
    \end{aligned}
\end{align*}
Therefore by \eqref{punt-G-conv} we have
\begin{equation}\label{v*}
\begin{cases}
v^*(s)\in \hat{V}_0,\\
-\div(\A ev^*(s))=\hat{h}(s).
\end{cases}
\end{equation}
Since, by \eqref{weak-fomulation-ellip-lap}, $\hat{v}_0(s)$ is a solution to \eqref{v*}, by uniqueness we have $v^*(s)=\hat{v}_0(s)$. Moreover, since the limit does not depend on the subsequence, the whole sequence satisfies 
\begin{equation}\label{weak-conv-hv-whole}
    \hv(s)\xrightharpoonup[\varepsilon\to 0^+]{}\hat{v}_0(s)\qquad\text{weakly in $\hat{V}_0$ for every $s\in \C_+$}.
\end{equation}

To prove the strong convergence we use $\hv(s)$ and $\hat{v}_0(s)$ as test function in \eqref{weak-fomulation-lap} and \eqref{weak-fomulation-ellip-lap}, respectively. By subtracting the two equalities, we obtain
\begin{equation*}
    (\A e\hv(s),e\hv(s))-(\A e\hat{v}_0(s),e\hat{v}_0(s))=(\hh(s),\hv(s))-(\hat{h}(s),\hat{v}_0(s))-\varepsilon^2s^2\|\hv(s)\|^2-\frac{\beta\varepsilon s}{\beta\varepsilon s+1}(\B e\hv(s),e\hv(s)),
\end{equation*}
from which we deduce
\begin{equation*}
    |(\A e\hv(s),e\hv(s))-(\A e\hat{v}_0(s),e\hat{v}_0(s))|\leq |(\hh(s),\hv(s))-(\hat{h}(s),\hat{v}_0(s))|+\varepsilon^2|s|^2\|\hv(s)\|^2+\beta\varepsilon |s|C_{\B}\|e\hv(s)\|^2.
\end{equation*}
By using again \eqref{punt-G-conv}, \eqref{bound-hv}, and \eqref{weak-conv-hv-whole}, we can deduce
\begin{equation}\label{ehv-conv}
    \lim_{\varepsilon\to 0^+}(\A e\hv(s),e\hv(s))=(\A e\hat{v}_0(s),e\hat{v}_0(s)) \quad\text{for every $s\in\C_+$}.
\end{equation}
Thanks to the coerciveness assumption \eqref{CB2Q}, the conclusion follows from the weak convergence \eqref{weak-conv-hv-whole} together with \eqref{ehv-conv}.
\end{proof}

Now we are in a position to prove the following result about the convergence in the space $L^2$ on the lines $\{s_1+is_2:s_2\in \R\}$.
\begin{proposition}\label{prop:conv-hv}
The functions $\hv$ and $\hat{v}_0$ satisfy \eqref{dadim}.
\end{proposition}

\begin{proof}
For every $s\in \C_+$, by using $\hv(s)$ as test function in \eqref{weak-fomulation-lap} we obtain
\begin{equation}\label{first-form}
    \frac{1}{\beta \varepsilon s+1}\Big(\beta \varepsilon^3s^3+\varepsilon^2 s^2+\beta \frac{((\A+\B)e\hv(s),e\hv(s))}{\|\hv(s)\|^2}\varepsilon s+\frac{(\A e\hv(s),e\hv(s))}{\|\hv(s)\|^2}\Big)\|\hv(s)\|^2=(\hh(s),\hv(s)).
\end{equation}
As before, we set
\begin{equation}\label{def-a}
    a:=\frac{(\A e\hv(s),e\hv(s))}{\|\hv(s)\|^2}\quad \text{and}\quad b:=\frac{((\A+\B)e\hv(s),e\hv(s))}{\|\hv(s)\|^2},
\end{equation}
and we observe that \eqref{coercQ} holds. Therefore, thanks to \eqref{dis-co}, Lemma \ref{comp-bound}, and \eqref{b-little} we can deduce
\begin{align*}
&\left|\frac{\beta \varepsilon^3s^3+\varepsilon^2 s^2+\beta b \varepsilon s+a}{\beta \varepsilon s+1}\right|\geq \left|\frac{\beta(\varepsilon s-z_0)}{\beta \varepsilon s+1}\right|\alpha\sqrt{\frac{b}{2}}\geq \beta\left|z_0\right|\alpha\sqrt{\frac{a }{2}}\geq \frac{\beta\alpha^2}{\sqrt{2}}\sqrt{a} & \quad\text{for $b>\frac{2}{3\beta^2}$},\\
&\left|\frac{\beta \varepsilon^3s^3+\varepsilon^2 s^2+\beta b \varepsilon s+a}{\beta \varepsilon s+1}\right|\geq \min\{\gamma,1\}a\geq \min\{\gamma,1\}\sqrt{a_0}\sqrt{a} &\quad\text{for $b\leq\frac{2}{3\beta^2}$},
\end{align*}
where in the first line we used the inequality $|z_0(\beta\varepsilon s+1)|\leq |\varepsilon s-z_0| $ for every $s\in \C_+$, which follows from the condition $z_0<0$.

As a consequence of these inequalities and of \eqref{first-form} there exists a positive constant $C=C(\alpha,\beta,\gamma,a_0)$ such that
\begin{align*}
   \|\hv(s)\|^2&=\left| \frac{\beta \varepsilon s+1}{\beta \varepsilon^3s^3+\varepsilon^2 s^2+\beta b\varepsilon s+a}\right||(\hh(s),\hv(s))|\leq \frac{C}{\sqrt{a}}\|\hh(s)\|\|\hv(s)\|\qquad\text{for every $s\in\C_+$}.
\end{align*}
Therefore, by using \eqref{def-a} and the coerciveness assumption \eqref{CB2Q}, we can write
\begin{equation*}
    \sqrt{c_{\A}}\|e\hv(s)\|\leq\sqrt{(\A e\hv(s),e\hv(s))}\leq C\|\hh(s)\|,
\end{equation*}
from which, recalling \eqref{KP-ineq}, we deduce
\begin{equation}\label{bound-hv-VD}
    \|\hv(s)\|_{\hat{V}_0}\leq(C_P+1) \frac{C}{\sqrt{c_{\A}}}\|\hh(s)\|\qquad\text{for every $s\in\C_+$}.
\end{equation}

By extending the function $\he$ to $(-\infty,0)$ with value $0$, we can write
\begin{align*}
    \hh(s)=\int_0^{+\infty}\e^{-st}\he(t)\de t=\int_{-\infty}^{+\infty}\e^{-st}\he(t)\de t=\mathcal{F}_t(\e^{-s_1t}\he(t))(s_2).
\end{align*}
Since for every $s=s_1+is_2\in\C_+$ the function  $t\mapsto\e^{-s_1t}\he(t)$ belongs to $L^2(\R;H)$, by the properties of the Fourier Transform we deduce that $s_2\mapsto\hh(s_1+is_2)$ belongs to $ L^2(\R;\hat{H})$ for every $\varepsilon>0$. Moreover, by using \eqref{conv-fe} and the Plancherel Theorem, we can write
\begin{align}\label{G-conv}
   \int_{-\infty}^{+\infty}\|\hh(s_1+is_2)&-\hat{h}(s_1+is_2)\|^2\de s_2=\int_{-\infty}^{+\infty}\|\mathcal{F}_t(\e^{-s_1t}(\he(t)-h(t)))(s_2)\|^2\de s_2\nonumber\\
    &=\int_{-\infty}^{+\infty}\|\e^{-s_1t}(\he(t)-h(t))\|^2\,\de t\leq \int_0^T\|\he(t)-h(t)\|^2\,\de t\xrightarrow[\varepsilon\to 0^+]{}0.
\end{align}
Since $\hv(s)\rightarrow \hat{v}_0(s)$ strongly in $\hat{V}_0$ by Lemma \ref{lem:conv-punt-hv} and $\hh(s)\rightarrow \hat{h}(s)$ strongly in $\hat{H}$ by \eqref{punt-G-conv}, thanks to \eqref{bound-hv-VD} and \eqref{G-conv} we can apply the Generalized Dominated Convergence Theorem to get the conclusion.
\end{proof}

\section{\texorpdfstring{$L^2$}{L2} convergence}
In this section we shall prove \eqref{conv-L2-f+g} and \eqref{conv-L2d-f+g} under the assumptions of Theorems \ref{main-thm-f+g} and \ref{main-f}. We begin by proving the following partial result.

\begin{proposition}\label{new-main-th}
Let $\{\he\}_{\varepsilon}\subset L^2(0,T;H)$ and $h\in L^2(0,T;H)$ be such that \eqref{conv-fe} holds. Let $\v$ and $v_0$ be the solutions to problems \eqref{weak-form-v} and \eqref{quas}. Then
\begin{equation*}
\v\xrightarrow[\varepsilon\to 0^+]{}v_0\qquad\text{strongly in $L^2(0,T;V)$}.
\end{equation*}
\end{proposition}

\begin{proof}
By the Plancherel Theorem we deduce from \eqref{inv-lap-four} and Proposition \ref{prop:conv-hv} that for every $s_1>0$ there exists a positive constant $C=C(s_1,T)$ such that
\begin{align*}
    \|\v -v_0\|^2_{L^2(0,T;V)}&=\int_0^T\|\v (t)-v_0(t)\|^2_{V}\,\de t=\int_0^T\|\mathcal{L}^{-1}(\hv-\hat{v}_0)(t)\|^2_{V}\,\de t\nonumber\\
    &\leq C(s_1,T)\int_{-\infty}^{+\infty}\|\mathcal{F}^{-1}_{s_2}(\hv(s_1+is_2)-\hat{v}_0(s_1+is_2))(t)\|^2_{V}\de t\nonumber\\
    & =C(s_1,T)\int_{-\infty}^{+\infty}\|\hv(s_1+is_2)-\hat{v}_0(s_1+is_2)\|^2_{\hat{V}_0}\de s_2\xrightarrow[\varepsilon\to 0^+]{}0,
\end{align*}
which concludes the proof.
\end{proof}

\begin{theorem}\label{uconvf+g}
Let us assume (H1), (H2), and \eqref{conv-dati-u}. Let $\u$ be the solution to the viscoelastic dynamic system \eqref{spaz2}, with
$\phie=\fe$ and $\gae=\gef$, and let $u_0$ be the solution to the stationary problem \eqref{def-ellip}. Then \eqref{conv-L2-f+g} holds.
\end{theorem}

\begin{proof}
Thanks to Lemma \ref{lemma-z2} it is enough to prove the theorem in the case $z=0$ and $\ze=0$ for every $\varepsilon>0$. 
We divide the proof into two steps.

\medskip
\textit{Step 1. The case $\ueu=0$}. We reduce the problem to the case of homogeneous initial conditions by considering the functions
\begin{align}
    \v(t):=\u(t)-\uez\qquad\text{and}\qquad v_0(t):=u_0(t)-u^0 \qquad \text{for a.e.\  $t\in[0,T]$}.\label{ve}
\end{align}
Let us define
\begin{alignat}{2}
    &q_{\varepsilon}(t):=\gef(t)+\div(\A e\uez)+\e^{-\frac{t}{\beta\varepsilon}}\div(\B e\uez)&&\qquad\text{for every $t\in[0,T]$},\label{def-ge}\\
    &q(t):=g(t)+\div(\A e u^0)&&\qquad\text{for every $t\in[0,T]$}\label{def-g0}.
\end{alignat}
Since $\ueu=0$, it is easy to see that $\v$ satisfies \eqref{weaksol2} with $\he=\fe$, $\lee=q_{\varepsilon}$, $\vez=0$, and $\veu=0$, while $v_0$ satisfies \eqref{weaksol2-ellip} with $h=f$ and $\ell=q$. By \eqref{conv-dati-u} and \eqref{ve}, to prove \eqref{conv-L2-f+g}
it is enough to show that
\begin{equation}
    \v\xrightarrow[\varepsilon\to 0^+]{}v_0\qquad\text{strongly in $L^2(0,T;V)$.}\label{conv-v-L2-f+g}
\end{equation}

In order to apply Proposition \ref{new-main-th}, we approximate the forcing terms of the problems for $\v$ and $v_0$ by means of functions in $H^1(0,T;H)$ and we consider the corresponding solutions $\tve$ and $\tilde{v}_0$, for which 
Proposition \ref{new-main-th} yields $\tve\to\tilde{v}_0$ strongly in $L^2(0,T;V)$ as $\varepsilon\to 0^+$. Finally we show that $\|\tve-\v\|_{L^2(0,T;V)}$ and $\|\tilde{v}_0-v_0\|_{L^2(0,T;V)}$ are small uniformly with respect to $\varepsilon$, and this leads to the proof of \eqref{conv-v-L2-f+g}.

Let us fix $\delta>0$. Thanks to the density of $H$ in $V'_0$ and to Lemma \ref{lem:appro} we can find $\psi\in H^1(0,T;H)$ and $h_{\A}^0,h_{\B}^0\in H$  such that
\begin{align}
   \|\psi-g\|_{W^{1,1}(0,T;V'_0)}<\delta ,\quad\|h_{\A}^0-\div(\A eu^0)\|_{V'_0}<\delta,\quad\|h_{\B}^0-\div(\B eu^0)\|_{V'_0}<\delta.\label{f*}
\end{align}
Thanks to (H1) and \eqref{conv-dati-u} there exist $\varepsilon_0=\varepsilon_0(\delta)\in(0,\tfrac{1}{\beta})$ such that 
\begin{alignat}{4}
     &\|\psi-\gef\|_{W^{1,1}(0,T;V'_0)}<\delta,&&\quad\|h_{\A}^0-\div(\A e\uez)\|_{V'_0}<\delta,&&\quad\|h_{\B}^0-\div(\B e\uez)\|_{V'_0}<\delta&&\quad\text{for every $\varepsilon\in(0,\varepsilon_0)$}.\label{f}
\end{alignat}

Let $\varphi_{\varepsilon},\varphi\colon [0,T]\rightarrow H$ be defined by
\begin{align}\label{def-fe}
\varphi_{\varepsilon}(t):=\psi(t)+h_{\A}^0+\e^{-\frac{t}{\beta\varepsilon}}h_{\B}^0\quad\text{and}\quad \varphi(t):=\psi(t)+h_{\A}^0\quad\text{for every $t\in[0,T]$}.
\end{align}
By \eqref{def-ge}, \eqref{def-g0}, \eqref{f*}, \eqref{f}, and \eqref{def-fe} for every $\varepsilon\in(0,\varepsilon_0)$ we obtain
\begin{align}
\|\varphi_{\varepsilon}-q_{\varepsilon}\|_{W^{1,1}(0,T;V'_0)}&\leq \|\psi-\gef\|_{W^{1,1}(0,T;V'_0)}+T\|h_{\A}^0-\div(\A e\uez)\|_{V'_0}\nonumber\\
    &\qquad+(\beta\varepsilon+1)\|h_{\B}^0-\div(\B e\uez)\|_{V'_0}\le (3+T)\delta,
    \label{5*}\\
    \|\varphi-q\|_{L^\infty(0,T;V'_0)}
&\leq (1+\tfrac{1}{T})\|\psi-g\|_{W^{1,1}(0,T;V'_0)}+\|h_{\A}^0-\div(\A e\uez)\|_{V'_0}\le (2+\tfrac{1}{T})\delta.
    \label{5**}
    \end{align}
Since $t\mapsto \e^{-\frac{t}{\beta\varepsilon}}\psi^0_{\B}$ converges to $0$ strongly in $L^2(0,T;H)$ as $\varepsilon\to 0^+$, by \eqref{def-fe} we have
\begin{equation}\label{convh}
     \varphi_{\varepsilon}\xrightarrow[\varepsilon\to 0^+]{} \varphi\qquad\text{strongly in $L^2(0,T;H)$}.
\end{equation}

Let $\tve$ be the solution to \eqref{weak-form-v} with $\he=\fe+ \varphi_{\varepsilon}$ and let $\tilde{v}_0$ be the solution to \eqref{quas} with $h=f+ \varphi$. By (H1) and \eqref{convh} we have
\begin{equation*}
     \fe+\varphi_{\varepsilon}\xrightarrow[\varepsilon\to 0^+]{}f+\varphi\qquad\text{strongly in $L^2(0,T;H)$},
\end{equation*}
hence Proposition \ref{new-main-th} yields
\begin{equation}\label{conv-tilde}
    \tve\xrightarrow[\varepsilon\to 0^+]{}\tilde{v}_0\qquad \text{strongly in $L^2(0,T;V)$}.
\end{equation}

To estimate the difference $\tve-\v$ we observe that it solves \eqref{weaksol2} with $\he=0$, $\lee=\varphi_{\varepsilon}-q_{\varepsilon}$, $\vez=0$, and $\veu=0$. Therefore, by Lemma \ref{lem:en-est} we have
\begin{equation}\label{4}
    \|\tve-\v\|_{L^2(0,T;V)}\leq \sqrt{C_E T}\|\varphi_{\varepsilon}-q_{\varepsilon}\|_{W^{1,1}(0,T;V'_0)}.
\end{equation}
To estimate the difference $\tilde{v}_0-v_0$
we observe that it solves \eqref{weaksol2-ellip} with $h=0$ and $\ell=\varphi-q$. Therefore by the Lax-Milgram Lemma we obtain
\begin{equation}\label{int-1}
    \|\tilde{v}_0-v_0\|_{L^{2}(0,T;V)}\leq  \tfrac{\sqrt{T}(C^2_P+1)}{c_{\A}}\|\varphi-q\|_{L^{\infty}(0,T;V'_0)}.
\end{equation}

By \eqref{5*}, \eqref{5**}, \eqref{4}, and \eqref{int-1} there exists a positive constant $C=C(\A,\B,\O,T)$ such that
\begin{equation*}
    \|\tve-\v\|_{L^2(0,T;V)}+ \|\tilde{v}_0-v_0\|_{L^{2}(0,T;V)}\leq C\delta,
\end{equation*}
hence
\begin{align*}
    \|\v-v_0\|_{L^2(0,T;V)}&\leq \|\v-\tve\|_{L^2(0,T;V)}+\|\tve-\tilde{v}_0\|_{L^2(0,T;V)}+\|\tilde{v}_0-v_0\|_{L^2(0,T;V)}\leq\|\tve-\tilde{v}_0\|_{L^2(0,T;V)}+C\delta.
\end{align*}
This inequality, together with \eqref{conv-tilde}, gives
\begin{equation*}
    \limsup_{\varepsilon\to 0^+} \|\v-v_0\|_{L^2(0,T;V)}\leq C\delta.
\end{equation*}
By the arbitrariness of $\delta>0$ we obtain
\eqref{conv-v-L2-f+g}, which concludes the proof of Step~1.

\medskip
\textit{Step 2. The general case}.
Let $\tue$ be the solution to \eqref{weaksol2} with $\he=\fe$, $\lee=\gef$, $\vez=\uez$, and $\veu=0$. By Step 1
\begin{equation}
    \tue\xrightarrow[\varepsilon\to 0^+]{}u_0\qquad\text{strongly in $L^2(0,T;V)$.}\label{conv-v-L2-f+g-tilde}
\end{equation}
The function $\u-\tue$ is the solution to \eqref{weaksol2} with all data equal to $0$ except $\veu$, which is now equal to $\ueu$. Therefore, Lemma \ref{lem:en-est} and \eqref{conv-dati-u} yield
\begin{equation*}
    \|\u-\tue\|_{L^{\infty}(0,T;V)}\leq \varepsilon\sqrt{C_E}\|\ueu\|\xrightarrow[\varepsilon\to 0^+]{}0,
\end{equation*}
which, together with \eqref{conv-v-L2-f+g-tilde}, gives \eqref{conv-L2-f+g}.
\end{proof}

In the following theorem, under the assumptions of Theorem \ref{main-f} we deduce \eqref{conv-L2d-f+g} from \eqref{conv-L2-f+g}.

\begin{theorem}\label{uconvf+gder}
Let us assume (H1), (H2), and \eqref{conv-dati-u}. Let $\u$ be the solution to the viscoelastic dynamic system \eqref{spaz2}, with $\phie=\fe$ and $\gae=\gef$, and let $u_0$ be the solution to the stationary problem \eqref{def-ellip}. Then \eqref{conv-L2d-f+g} holds.
\end{theorem}

\begin{proof}
Thanks to Lemma \ref{lemma-z2} we can suppose $z=0$ and $\ze=0$ for every $\varepsilon>0$. It is convenient to extend the data of our problem to the interval $[0,2T]$ by setting
\begin{equation*}
    \fe(t):=0,\quad f(t):=0,\quad \gef(t):=\gef(T),\quad g(t):=g(T)\qquad \text{for every $t\in (T,2T]$}.
\end{equation*}
Since (H1) holds, it is clear that $\{\fe\}_{\varepsilon}\subset L^2(0,2T;H)$, $\{\gef\}_{\varepsilon}\in H^1(0,2T;V'_0)$,
\begin{equation}\label{newconvfor}
    \fe\xrightarrow[\varepsilon\to 0^+]{}f\quad \text{strongly in $L^2(0,2T;H)$}\quad \text{and} \quad \gef\xrightarrow[\varepsilon\to 0^+]{}g\quad \text{strongly in $W^{1,1}(0,2T;V'_0)$.}
\end{equation}
Moreover, the solution to \eqref{spaz2} on $[0,2T]$ with
the extended data is an extension of $\u$, which is still denoted by $\u$. Similarly, the solution to \eqref{def-ellip} on $[0,2T]$ is still denoted by $u_0$. Since \eqref{newconvfor} holds, Theorem \ref{uconvf+g} gives
\begin{equation}\label{tildeconv}
    \u\xrightarrow[\varepsilon\to 0^+]{}u_0\quad \text{strongly in $L^2(0,2T;V)$}.
\end{equation}
We further extend $\u$ to $\R$ by setting $\u(t)=0$ for every $t\in \R\setminus[0,2T]$, and we define
\begin{equation*}
\w(t):=\int_0^t\frac{1}{\beta\varepsilon}\e^{-\frac{t-\tau}{\beta\varepsilon}}e\u(\tau)\de\tau=(\rho_{\varepsilon}*e\u)(t)\quad \text{for every $t\in\R$},
\end{equation*}
where $\rho_{\varepsilon}$ is as in \eqref{rho}. By  the properties of convolutions and \eqref{tildeconv} we get
\begin{equation}\label{eww}
    e\u-\w\xrightarrow[\varepsilon\to 0^+]{} 0\quad \text{strongly in $L^2(\R;\tilde{H})$}.
\end{equation}

Thanks to \eqref{tildeconv} and \eqref{eww}, by using \eqref{spaz2} and \eqref{def-ellip} we obtain
\begin{equation}\label{convusec}
    \varepsilon^2\ddu\xrightarrow[\varepsilon\to 0^+]{}0\quad \text{strongly in $L^2(0,2T;V'_0)$}.
\end{equation}
Since
\begin{equation*}
    \varepsilon^2\du(t)=\varepsilon^2\ueu+\varepsilon^2\int_0^t\ddu(\tau)\de\tau\quad\text{for every $t\in[0,2T]$},
\end{equation*}
\eqref{conv-dati-u} and \eqref{convusec} imply
\begin{equation}\label{convufis}
    \varepsilon^2\du\xrightarrow[\varepsilon\to 0^+]{}0\quad \text{strongly in $L^2(0,2T;V'_0)$}.
\end{equation}
By \eqref{tildeconv} and \eqref{convufis} there exists a sequence $\varepsilon_j\xrightarrow[]{}0^+$ such that for a.e.\ $t\in[0,2T]$ we have
\begin{alignat}{2}
   \uk(t)&\xrightarrow[j\to +\infty]{} u_0(t)&&\qquad \text{strongly in $V$},\label{jkconv}\\
   \varepsilon^2_j\duk(t)&\xrightarrow[j\to +\infty]{} 0&&\qquad \text{strongly in $V'_0$}.\label{jkconv2}
\end{alignat}
We choose $T_0\in (T,2T)$ such that \eqref{jkconv} and \eqref{jkconv2} hold at $t=T_0$. This implies
\begin{equation}\label{t0}
    \varepsilon_{j}^2(\duk(T_0),\uk(T_0))=\langle\varepsilon_{j}^2\duk(T_0),\uk(T_0)\rangle\xrightarrow[j\to +\infty]{}0.
\end{equation}

Since $\ze=0$ for a.e.\ $t\in[0,T_0]$ we can use $\u(t)\in V_0$ as test function in \eqref{spaz2}. Then we integrate by parts in time on the interval $(0,T_0)$ to obtain
\begin{align*}
-\varepsilon^2_{j}&\int_0^{T_0}\|\duk(t)\|^2\de t+\int_0^{T_0}(\A e\uk(t),e\uk(t))\de t+\int_0^{T_0}(\B(e\uk(t)-w_{\varepsilon_{j}}(t)),e \uk(t))\de t\nonumber\\
&=\int_0^{T_0}(\fuk(t),\uk(t))\de t+\int_0^{T_0}\langle  g_{\varepsilon_j}(t),\uk(t)\rangle\de t-\varepsilon^2_j(\duk(T_0),\uk(T_0))+\varepsilon^2_j(u^0_{\varepsilon_j},u^1_{\varepsilon_j}).
\end{align*}
Thanks to \eqref{def-ellip}, \eqref{conv-dati-u}, \eqref{newconvfor}, \eqref{tildeconv}, \eqref{eww}, and \eqref{t0} the first term on the left-hand side of the previous equation tends to $0$ as $j\to +\infty$. Since $T_0>T$ we have
\begin{equation*}
    \varepsilon^2_{j}\int_0^{T}\|\dot{u}_{\varepsilon_{j}}(t)\|^2\de t\xrightarrow[j\to+\infty]{}0.
\end{equation*}

By the arbitrariness of the sequence $\{\varepsilon_j\}_j$ we have
\begin{equation*}
    \varepsilon^2\int_0^{T}\|\du(t)\|^2\de t\xrightarrow[\varepsilon\to 0^+]{}0,
\end{equation*}
which concludes the proof.
\end{proof}

We now use Theorems \ref{uconvf+g} and \ref{uconvf+gder} to obtain \eqref{conv-L2-f+g} and \eqref{conv-L2d-f+g} under the assumptions of Theorem \ref{main-thm-f+g}.

\begin{theorem}\label{main-thm-f+g-part1}
Let us assume (H1)--(H3). Let $\u$ be the solution to the viscoelastic dynamic system \eqref{spaz} and let $u_0$ be the solution to the stationary problem \eqref{def-ellip}. Then \eqref{conv-L2-f+g} and \eqref{conv-L2d-f+g} hold.
\end{theorem}

\begin{proof}
Thanks to Lemma \ref{lemma-z} we can suppose $z=0$ and $\ze=0$ for every $\varepsilon>0$.
Let $p_{\varepsilon}$ be defined by \eqref{g0}.
Since $\ze=0$, by Remark \ref{corris} the function $\u$ solves \eqref{weaksol2} with $\he=\fe$, $\lee=\gef-p_{\varepsilon}$, $\vez=\uein(0)$, and 
$\veu=\duein(0)$.
To obtain \eqref{conv-L2-f+g} and \eqref{conv-L2d-f+g} we cannot apply
Theorems \ref{uconvf+g} and \ref{uconvf+gder} directly, because 
$\{p_{\varepsilon}\}_{\varepsilon}$ does not converge to
$0$ in $W^{1,1}(0,T;V'_0)$ as $\varepsilon\to 0^+$ and, in general, $p_{\varepsilon}\notin L^2(0,T;H)$.

To overcome this difficulty we construct a family
$\{q_\varepsilon\}_{\varepsilon}\subset H^1(0,T;H)$ such that $\|{q_\varepsilon-p_\varepsilon}\|_{W^{1,1}(0,T;V'_0)}$ is uniformly small and $q_\varepsilon\to 0$ strongly in $L^2(0,T;H)$ as $\varepsilon\to 0^+$. 
Then we can apply Theorems \ref{uconvf+g} and \ref{uconvf+gder} to the solutions $\v$ to \eqref{weaksol2} with $p_\varepsilon$ replaced by $q_\varepsilon$, obtaining
that $\v\to u_0$ strongly in $L^2(0,T;V)$ and $\varepsilon\dv\to 0$ strongly in $L^2(0,T;H)$. Finally, we show that $\|\v-\u\|_{L^2(0,T;V)}$ and $\varepsilon\|\dv-\du\|_{L^2(0,T;H)}$ are small uniformly with respect to $\varepsilon$, and this leads to the proof of \eqref{conv-L2-f+g} and \eqref{conv-L2d-f+g}.

To construct $q_\varepsilon$ we
consider $g^0_{\varepsilon}$ introduced in \eqref{g0} and we define 
\begin{equation*}
    \tilde{g}^0_{\varepsilon}:=\int_{-\infty}^0\frac{1}{\beta\varepsilon}\e^{\frac{\tau}{\beta\varepsilon}}\div(\B e\uin(\tau))\de\tau=(\rho_{\varepsilon}*\div(\B e\uin))(0),
\end{equation*}
where $\rho_{\varepsilon}$ is the convolution kernel in \eqref{rho}. By (H3) we have $\div(\B e\uin)\in C^0((-\infty,0];V'_0)$, hence the properties of convolutions imply
\begin{equation}\label{limitatezza-inf}
  \tilde{g}^0_{\varepsilon}\xrightarrow[\varepsilon\to 0^+]{}g^0:=\div(\B e \uin(0))\quad\text{strongly in $V'_0$.}
\end{equation}

Since
\begin{align*}
    \|g^0_{\varepsilon}-\tilde{g}^0_{\varepsilon}\|_{V'_0}&\leq \int_{-\infty}^{-a}\frac{1}{\beta\varepsilon}\e^{\frac{\tau}{\beta\varepsilon}}\|\div(\B(e\uein(\tau))\|_{V'_0}\de\tau+\int_{-\infty}^{-a}\frac{1}{\beta\varepsilon}\e^{\frac{\tau}{\beta\varepsilon}}\|\div(\B(e\uin(\tau))\|_{V'_0}\de\tau\\
    &\hspace{2cm}+\|\div(\B (e\uein-e\uin))\|_{L^{\infty}(-a,0;V'_0)},
\end{align*}
thanks to (H3) we have $g^0_{\varepsilon}-\tilde{g}^0_{\varepsilon}\rightarrow 0$ strongly in $V'_0$ as $\varepsilon\to 0^+$, hence \eqref{limitatezza-inf} implies 
\begin{equation}\label{limitatezza*-inf}
 g^0_{\varepsilon}\xrightarrow[\varepsilon\to 0^+]{}g^0\quad\text{strongly in $V'_0$.}
\end{equation}

Let us fix $\delta>0$. By the density of $H$ in $V'_0$ we can find $h^0\in H$  such that
$\|h^0-g^0\|_{V'_0}<\delta$.
By \eqref{limitatezza*-inf} there exists $\varepsilon_0=\varepsilon_0(\delta)\in(0,\tfrac{1}{\beta})$ such that 
\begin{equation}
     \|h^0-g^0_{\varepsilon}\|_{V'_0}<\delta\quad\text{for every $\varepsilon\in(0,\varepsilon_0)$}.\label{f-1f}
\end{equation}
Let $q_{\varepsilon}\in H^1(0,T;H)$ be defined by
$q_{\varepsilon}(t):=\e^{-\frac{t}{\beta\varepsilon}}h^0$
for every $t\in[0,T]$.
Then
\begin{equation}\label{convhf}
     q_{\varepsilon}\xrightarrow[\varepsilon\to 0^+]{} 0\qquad\text{strongly in $L^2(0,T;H)$}.
\end{equation}
Since $p_\varepsilon(t)=\e^{-\frac{t}{\beta\varepsilon}}g^0_\varepsilon$, by \eqref{f-1f} we have also
\begin{equation}\label{convhf*}
     \|q_{\varepsilon}-p_{\varepsilon}\|_{W^{1,1}(0,T;V'_0)}\le (\beta\varepsilon+1)
     \|h^0-g^0_{\varepsilon}\|_{V'_0}
     \le 2\delta\quad\text{for every }\varepsilon\in(0,\varepsilon_0).
\end{equation}

Let $\v$ be the solution to
\eqref{weaksol2} with $\he=\fe-q_{\varepsilon}$, $\lee=\gef$, $\vez=\uein(0)$, and 
$\veu=\duein(0)$.
By (H1) and \eqref{convhf} we have
\begin{equation*}
   \fe-q_{\varepsilon}\xrightarrow[\varepsilon\to 0^+]{} f\quad\text{strongly in $L^2(0,T;H)$} \quad\text{and}\quad \gef\xrightarrow[\varepsilon\to 0^+]{}g\quad\text{strongly in $W^{1,1}(0,T;V'_0)$}.
\end{equation*}
By (H3) we have 
\begin{equation*}
   \uein(0)\xrightarrow[\varepsilon\to 0^+]{} \uin(0)\quad\text{strongly in $V$} \quad\text{and}\quad \varepsilon\duein(0)\xrightarrow[\varepsilon\to 0^+]{}0\quad\text{strongly in $H$}.
\end{equation*}
Therefore we can apply Theorems \ref{uconvf+g} and \ref{uconvf+gder} to obtain
\begin{equation}\label{conv-tildef}
\v\xrightarrow[\varepsilon\to 0^+]{}u_0\quad \text{strongly in $L^2(0,T;V)$}\quad\text{and}\quad
    \varepsilon\dv\xrightarrow[\varepsilon\to 0^+]{}0\quad \text{strongly in $L^2(0,T;H)$}.
\end{equation}

To estimate the difference $\v-\u$ we observe that it solves \eqref{weaksol2} with $\he=0$, $\lee=p_{\varepsilon}-q_\varepsilon$,
$\vez=0$, and 
$\veu=0$.
Therefore, by Lemma \ref{lem:en-est} and \eqref{convhf*} we have
\begin{equation}\label{4f}
    \varepsilon^2\|\dv-\du\|^2_{L^2(0,T;H)}+\|\v-\u\|^2_{L^2(0,T;V)}\leq C_E \|q_{\varepsilon}-p_{\varepsilon}\|^2_{W^{1,1}(0,T;V'_0)}
    \leq 4C_E \delta^2.
\end{equation}
Since by \eqref{4f}
\begin{align*}
    \|\u-u_0\|_{L^2(0,T;V)}&\leq \|\u-\v\|_{L^2(0,T;V)}+\|\v-u_0\|_{L^2(0,T;V)}
    \leq\|\v-u_0\|_{L^2(0,T;V)}
    +2\sqrt{C_E}\delta,\\
    \varepsilon\|\du\|_{L^2(0,T;H)}&\leq \varepsilon\|\du-\dv\|_{L^2(0,T;H)}+\varepsilon\|\dv\|_{L^2(0,T;H)}\leq\varepsilon\|\dv\|_{L^2(0,T;H)}+2\sqrt{C_E}\delta,
\end{align*}
thanks to \eqref{conv-tildef} we have
\begin{equation*}
    \limsup_{\varepsilon\to 0^+}\|\u-u_0\|_{L^2(0,T;V)}\leq 2\sqrt{C_E}\delta\quad\text{and}\quad  \limsup_{\varepsilon\to 0^+}\varepsilon\|\du\|_{L^2(0,T;H)}\leq 2\sqrt{C_E}\delta.
\end{equation*}
By the arbitrariness of $\delta>0$ we obtain \eqref{conv-L2-f+g} and \eqref{conv-L2d-f+g}, which concludes the proof.
\end{proof}

\section{The local uniform convergence}
In this section we shall prove \eqref{th-4-g-f+g} and \eqref{th-3-g-f+g} under the assumptions of Theorems \ref{main-thm-f+g} and \ref{main-f}. The proof is based on the following lemma. 

\begin{lemma}\label{regular-f+g}
Let $\{\lee\}_{\varepsilon}\subset H^{1}(0,T;V'_0)$ and $\ell\in W^{1,1}(0,T;V'_0)$ be such that
\begin{equation}
    \lee\xrightarrow[\varepsilon\to 0^+]{}\ell\quad \text{strongly in $W^{1,1}(\eta,T;V'_0)$ for every $\eta\in(0,T)$}.\label{conv-gamma*}
\end{equation}
Let $\v$ be a solution to the viscoelastic dynamic system \eqref{weaksol2} with $\he=0$ and arbitrary initial data. Moreover, let $v_0$ be the solution to the stationary problem \eqref{weaksol2-ellip} with $h=0$. We assume that
\begin{align}
    \v\xrightarrow[\varepsilon\to 0^+]{}v_0\qquad\text{strongly in $L^2(0,T;V)$,}\label{conv-med1-*}\\
     \varepsilon\dv\xrightarrow[\varepsilon\to 0^+]{}0\qquad\text{strongly in $L^2(0,T;H)$.}\label{conv-med2-*}
\end{align}
Then
\begin{align}
    \v\xrightarrow[\varepsilon\to 0^+]{}v_0\qquad\text{strongly in $L^{\infty}(\eta,T;V)$ for every $\eta\in(0,T)$,}\label{tt}\\
     \varepsilon\dv\xrightarrow[\varepsilon\to 0^+]{}0\qquad\text{strongly in $L^{\infty}(\eta,T;H)$ for every $\eta\in(0,T)$.}\label{tt2}
\end{align}
\end{lemma}

\begin{proof}
We divide the proof into two steps.

\medskip
{\it Step 1.} Let us assume $\lee=\ell\in H^2(0,T;V'_0)$ for every $\varepsilon>0$. By Lemma \ref{reg-ell-Hk} (with $z=0$) we have $v_0\in H^2(0,T;V)$, hence recalling \eqref{weaksol2-ellip} we get
\begin{align}\label{iper-u0**}
    \varepsilon^2\ddot v_0(t) &-\div((\A+\B) e v_0(t))+\int_0^t\frac{1}{\beta\varepsilon}\e^{-\frac{t-\tau}{\beta\varepsilon}}\div(\B e v_0(\tau))\de\tau\nonumber\\
    &= \varepsilon^2\ddot v_0(t)+\ell(t)-\div(\B e v_0(t))+\int_0^t\frac{1}{\beta\varepsilon}\e^{-\frac{t-\tau}{\beta\varepsilon}}\div(\B e v_0(\tau))\quad\text{for a.e.\ }t\in [0,T].
\end{align}
Now we define $\vb:=\v-v_0$ and observe that by \eqref{conv-med1-*} and \eqref{conv-med2-*} we have
\begin{alignat}{2}
\vb&\xrightarrow[\varepsilon\to 0^+]{}0&&\qquad\text{strongly in $L^2(0,T;V)$},\label{th-1-be-*}\\
\varepsilon\dvb&\xrightarrow[\varepsilon\to 0^+]{}0&&\qquad\text{strongly in $L^2(0,T;H)$}.\label{th-2-be-*} 
\end{alignat} 
Let us consider
 \begin{equation*}
   q_{\varepsilon}(t):=\div(\B ev_0(t))-\int_0^t \frac{1}{\beta\varepsilon}\e^{-\frac{t-\tau}{\beta\varepsilon}}\div(\B e v_0(\tau))\de \tau.
\end{equation*}
Since $\v$ satisfies \eqref{weaksol2} with $\he=0$, by \eqref{iper-u0**} the function $\vb$ satisfies \eqref{weaksol2} with $\he=-\varepsilon^2\ddot v_0$ and $\lee=q_{\varepsilon}$. After two integrations by parts in time we deduce
\begin{align*}
    \int_0^t\frac{1}{\beta\varepsilon} \e^{-\frac{t-\tau}{\beta\varepsilon}}\div(\B e v_0(\tau))\de \tau=\div(\B e v_0(t))&-\e^{-\frac{t}{\beta\varepsilon}}\div(\B e v_0(0))-\beta\varepsilon\div(\B e \dot{v}_0(t))\nonumber\\
    &+\beta\varepsilon\e^{-\frac{t}{\beta\varepsilon}}\div(\B e \dot{v}_0(0))+\beta\varepsilon\int_0^t\e^{-\frac{t-\tau}{\beta\varepsilon}}\div(\B e \ddot{v}_0(\tau))\de \tau,
\end{align*}
hence 
\begin{alignat}{2}
 q_{\varepsilon}&\xrightarrow[\varepsilon\to 0^+]{}0&&\qquad\text{strongly in $W^{1,1}(\eta,T;V'_0)$ for every $\eta\in(0,T)$}.\label{2-be*}
\end{alignat}

Now we fix $\delta\in(0,T)$, and we consider $\eta\in (0,\delta)$ and $\zeta\in (\eta,\delta)$. We define the family of functions $\{\wb\}_{\varepsilon}\subset H^1(0,T;\tilde{H})$ by
\begin{equation*}
    \wb(t):=\int_0^t\frac{1}{\beta\varepsilon}\e^{-\frac{t-\tau}{\beta\varepsilon}}e\vb(\tau)\de\tau=(\rho_{\varepsilon}*e\vb)(t)\quad\text{for every $t\in[0,T]$},
\end{equation*}
where $\rho_{\varepsilon}$ is defined by \eqref{rho} and $\vb$ is extended to $\R$ by setting $\vb(t)=0$ on $\R\setminus[0,T]$. By properties of convolutions we have
\begin{equation}\label{conv-w*}
    e\vb-\wb\xrightarrow[\varepsilon\to 0^+]{}0\quad\text{strongly in $L^2(0,T;\tilde{H})$}.
\end{equation}

By the energy-dissipation balance \eqref{energy} of Proposition \ref{diss-bal}, for every $t\in [\eta,T]$ and $s\in (\eta,\zeta)$ we can write
\begin{align}\label{eq-en-t0}
    \frac{\varepsilon^2}{2}&\|\dvb(t)\|^2+\frac{1}{2}(\A e\vb(t),e\vb(t))+\frac{1}{2}(\B( e\vb(t)-\wb(t),e\vb(t)-\wb(t))+\beta\varepsilon\int_s^{t}(\B\dwb(\tau),\dwb(\tau))\de\tau\nonumber\\
    &=\frac{\varepsilon^2}{2}\|\dvb(s)\|^2+\frac{1}{2}(\A e\vb(s),e\vb(s))+\frac{1}{2}(\B( e\vb(s)-\wb(s),e\vb(s)-\wb(s))+\W_{\varepsilon}(t,s),
\end{align}
where the work is defined by
\begin{align*}
  \W_{\varepsilon}(t,s)&=\langle q_{\varepsilon}(t),\vb(t)\rangle-\langle q_{\varepsilon}(s),\vb(s)\rangle-\int_s^{t} \langle\dot q_{\varepsilon}(\tau),\vb(\tau)\rangle\de \tau-\varepsilon\int_s^{t}(\ddot v_0(\tau),\varepsilon\dvb(\tau))\de \tau.
\end{align*}
Now we take the mean value with respect to $s$ of all terms of \eqref{eq-en-t0} on $(\eta,\zeta)$, and we pass to the supremum with respect to $t$ on $[\eta,T]$. Thanks to \eqref{KP-ineq} and \eqref{CB2Q} we deduce
\begin{align}\label{e*}
    \frac{\varepsilon^2}{2}\|\dvb\|^2_{L^{\infty}(\eta,T;H)}&+\frac{c_{\A}}{2(C^2_P+1)}\|\vb\|_{L^{\infty}(\eta,T;V)}^2\leq\frac{\varepsilon^2}{2}\dashint_{\eta}^{\zeta}\|\dvb(s)\|^2\de s+\frac{C_{\A}}{2}\dashint_{\eta}^{\zeta}\| \vb(s)\|^2_V\de s\nonumber\\
    &+\frac{C_{\B}}{2}\dashint_{\eta}^{\zeta}\|e\vb(s)-\wb(s)\|^2\de s+\dashint_{\eta}^{\zeta}\sup_{t\in[\eta,T]}|\W_{\varepsilon}(t,s)|\de s.
\end{align}
Notice that for every $s\in(\eta,\zeta)$ we have
\begin{align*}
    \sup_{t\in[\eta,T]}|\W_{\varepsilon}(t,s)|&\leq \big(2\|q_{\varepsilon}\|_{L^{\infty}(\eta,T;V'_0)}+\|\dot q_{\varepsilon}\|_{L^{1}(\eta,T;V'_0)}\big)\|\vb\|_{L^{\infty}(\eta,T;V)}+\varepsilon\|\ddot v_0\|_{L^1(\eta,T;H)}\|\varepsilon\dvb\|_{L^{\infty}(\eta,T;H)}\\
    &\leq \big(3+\tfrac{2}{T}\big)\|q_{\varepsilon}\|_{W^{1,1}(\eta,T;V'_0)}\|\vb\|_{L^{\infty}(\eta,T;V)}+\varepsilon\|\ddot v_0\|_{L^1(\eta,T;H)}\|\varepsilon\dvb\|_{L^{\infty}(\eta,T;H)},
\end{align*}
hence thanks to the Young Inequality and \eqref{e*} there exists a positive constant $C=C(\A,\B,\O,T)$ such that
\begin{align}\label{e**}
    \varepsilon^2\|\dvb\|^2_{L^{\infty}(\eta,T;H)}&+\|\vb\|_{L^{\infty}(\eta,T;V)}^2\leq C\Big(\varepsilon^2\dashint_{\eta}^{\zeta}\|\dvb(s)\|^2\de s+\dashint_{\eta}^{\zeta}\| \vb(s)\|^2_V\de s\nonumber\\
    &+\dashint_{\eta}^{\zeta}\|e\vb(s)-\wb(s)\|^2\de s+\|q_{\varepsilon}\|^2_{W^{1,1}(\eta,T;V'_0)}+\varepsilon^2\|\ddot v_0\|^2_{L^1(\eta,T;H)}\Big).
\end{align}

By passing to the limit in \eqref{e**} as $\varepsilon\to 0^+$, thanks to \eqref{th-1-be-*}, \eqref{th-2-be-*}, \eqref{2-be*}, and \eqref{conv-w*} we obtain
\begin{equation*}
   \varepsilon\|\dvb\|_{L^{\infty}(\eta,T;H)}+\|\vb\|_{L^{\infty}(\eta,T;V)}\xrightarrow[\varepsilon\to 0^+]{}0,
\end{equation*}
which, by the definition of $\vb$, concludes the proof of \eqref{tt} and \eqref{tt2} in the case $\ell\in H^2(0,T;V'_0)$.

\medskip
{\it Step 2.} In the general case $\ell\in W^{1,1}(0,T;V'_0)$ we use an approximation argument. Given $\delta>0$, by Lemma \ref{densita} there exists a function $\psi\in H^2(0,T;H)$ such that 
\begin{equation}\label{etaf}
\|\psi-\ell\|_{W^{1,1}(0,T;V'_0)}<\delta.    
\end{equation}
Thanks to \eqref{conv-gamma*} for every $\sigma\in(0,T)$ there exists a positive number $\varepsilon_0=\varepsilon_0(\delta,\sigma)$ such that
\begin{equation}\label{etas}
\|\psi-\lee\|_{W^{1,1}(\sigma,T;V'_0)}<\delta\qquad\text{for every $\varepsilon\in(0,\varepsilon_0)$}.    
\end{equation}
Let $\tve$ be the solution to \eqref{weaksol2} in the interval $[\sigma,T]$ with $\he=0$, $\lee=\psi$, $\tve(\sigma)=\v(\sigma)$, and $\dtve(\sigma)=\dv(\sigma)$, and let $\tilde{v}_0$ be the solution to \eqref{weaksol2-ellip} in the interval $[0,T]$ with $h=0$ and $\ell=\psi$. By applying Step 1 in the interval $[\sigma,T]$ we obtain
\begin{alignat}{2}
    \tve&\xrightarrow[\varepsilon\to 0^+]{}\tilde{v}_0&&\qquad\text{strongly in $L^{\infty}(\eta,T;V)$ for every $\eta\in(\sigma,T)$},\label{convft}\\
   \varepsilon\dtve&\xrightarrow[\varepsilon\to 0^+]{}0&&\qquad\text{strongly in $L^{\infty}(\eta,T;H)$ for every $\eta\in(\sigma,T)$}.\label{convdt}
\end{alignat}

We set $\bar{v}_0:=\tilde{v}_0-v_0$ and $\vb:=\tve-\v$. We observe that $\bar{v}_0$ is the solution to \eqref{weaksol2-ellip} with $h=0$ and $\ell$ replaced by $\psi-\ell$, hence by the Lax-Milgram Lemma we get
\begin{equation}\label{stima-u0-f}
    \|\bar{v}_0\|_{L^{\infty}(0,T;V)}\leq \tfrac{C^2_P+1}{c_{\A}} \|\psi-\ell\|_{L^{\infty}(0,T;V'_0)}\leq\tfrac{C^2_P+1}{c_{\A}}(1+\tfrac{1}{T}) \|\psi-\ell\|_{W^{1,1}(0,T;V'_0)}.
\end{equation}
Moreover, $\vb$ is the solution to \eqref{weaksol2} in the interval $[\sigma,T]$ with $\he=0$, $\lee$ replaced by $\psi-\lee$, and homogeneous initial conditions.  Thanks to Lemma \ref{lem:en-est} we obtain
\begin{equation}\label{en-esi-2}
    \varepsilon\|\dvb\|^2_{L^{\infty}(\sigma,T;H)}+\|\vb\|^2_{L^{\infty}(\sigma,T;V)}\leq C_E\|\psi-\lee\|^2_{W^{1,1}(\sigma,T;V'_0)}.
\end{equation}
By combining \eqref{etaf}, \eqref{etas}, \eqref{stima-u0-f}, and \eqref{en-esi-2}, we can find a positive constant $C=C(\A,\B,\O,T)$ such that
\begin{equation}\label{en-esi-2*}
    \varepsilon\|\dvb\|_{L^{\infty}(\sigma,T;H)}+\|\vb\|_{L^{\infty}(\sigma,T;V)}+\|\bar{v}_0\|_{L^{\infty}(\sigma,T;V)}\leq C\delta.
\end{equation}
Since for every $\eta\in (\sigma,T)$ we have 
\begin{align*}
    \|\v-v_0\|_{L^{\infty}(\eta,T;V)}&\leq \|\vb\|_{L^{\infty}(\eta,T;V)}+\|\tve-\tilde{v}_0\|_{L^{\infty}(\eta,T;V)}+\|\bar{v}_0\|_{L^{\infty}(\eta,T;V)},\\
    \varepsilon\|\dv\|_{L^{\infty}(\eta,T;H)}&\leq \varepsilon\|\dvb\|_{L^{\infty}(\eta,T;H)}+\varepsilon\|\dtve\|_{L^{\infty}(\eta,T;H)},
\end{align*}
thanks to \eqref{convft}, \eqref{convdt}, and \eqref{en-esi-2*} we obtain
\begin{equation*}
    \limsup_{\varepsilon\to 0^+}\|\v-v_0\|_{L^{\infty}(\eta,T;V)}\leq C\delta\quad\text{and}\quad \limsup_{\varepsilon\to 0^+}\|\varepsilon\dv\|_{L^{\infty}(\eta,T;H)}\leq C\delta,
\end{equation*}
for every $\eta\in(\sigma,T)$. By the arbitrariness of $\delta>0$ and $\sigma>0$ we conclude.
\end{proof}

Now we are in a position to prove \eqref{th-4-g-f+g} and \eqref{th-3-g-f+g}. 

\begin{theorem}\label{teo-eta-L2}
Let us assume (H1), (H2), \eqref{conv-dati-u}, and $\fe=0$ for every $\varepsilon>0$. Let $\u$ be the solution to the viscoelastic dynamic system \eqref{spaz2}, with $\phie=0$ and $\gae=\gef$, and let $u_0$ be the solution to the stationary problem \eqref{def-ellip}, with $f=0$. Then \eqref{th-4-g-f+g} and \eqref{th-3-g-f+g} hold.
\end{theorem}

\begin{proof}
By Theorems \ref{uconvf+g} and \ref{uconvf+gder} we obtain \eqref{conv-L2-f+g} and \eqref{conv-L2d-f+g}. Since $\fe=0$ and $\gef\rightarrow g$ strongly in $W^{1,1}(0,T;V'_0)$ as $\varepsilon\to 0^+$ by (H1), we can apply Lemma \ref{regular-f+g} to conclude.
\end{proof}

\begin{theorem}\label{teo-eta-inf}
Let us assume (H1)--(H3) and $\fe=0$ for every $\varepsilon>0$. Let $\u$ be the solution to the viscoelastic dynamic system \eqref{spaz} and let $u_0$ be the solution to the stationary problem \eqref{def-ellip}, with $f=0$. Then \eqref{th-4-g-f+g} and \eqref{th-3-g-f+g} hold.
\end{theorem}

\begin{proof}
Thanks to Lemma \ref{lemma-z} we can suppose $z=0$ and $\ze=0$ for every $\varepsilon>0$. By Theorem \ref{main-thm-f+g-part1} we obtain \eqref{conv-L2-f+g} and \eqref{conv-L2d-f+g}. Since $\u$ is a solution to \eqref{spaz} with $\fe=0$, by Remark \ref{corris} it solves \eqref{weaksol2} with $\he=0$ and $\lee=\gef-p_{\varepsilon}$, where $p_{\varepsilon}$ is defined by \eqref{g0}. Since
\begin{equation*}
    \gef-p_{\varepsilon}\xrightarrow[\varepsilon\to 0^+]{}g\quad\text{strongly in $W^{1,1}(\eta,T;V'_0)$ for every $\eta\in(0,T)$,} 
\end{equation*}
we can apply Lemma \ref{regular-f+g} to conclude.
\end{proof}

Finally we can prove Theorems \ref{main-thm-f+g} and \ref{main-f}.

\begin{proof}[Proof of Theorem \ref{main-thm-f+g}] It is enough to combine Theorems \ref{teo-compat}, \ref{main-thm-f+g-part1}, and \ref{teo-eta-inf}.
\end{proof}

\begin{proof}[Proof of Theorem \ref{main-f}] It is enough to combine Theorems \ref{uconvf+g}, \ref{uconvf+gder}, and \ref{teo-eta-L2}.
\end{proof}

\appendix
\section{}

Throughout this section we fix $a_0>0$, $b_0>0$, and $c_1\geq c_0>1$. For every $a,b$ with
\begin{alignat}{2}\label{ip1}
    c_0a\leq b\leq c_1  a,\qquad b\geq b_0,\qquad a\geq a_0,
\end{alignat}
we consider the polynomial $p(z):=\beta z^3+z^2+\beta b z+a$ depending on the complex variable $z$. The following result about the roots of this polynomial is used in the proof of Lemma \ref{lem:ex-lax} and Proposition \ref{prop:conv-hv}.

\begin{lemma}\label{comp-bound}
There exists a positive constant $\alpha=\alpha(\beta,a_0,b_0,c_0,c_1)$ such that, for every $a,b\in\R$ satisfying \eqref{ip1}, the roots of the polynomial $p$ have real parts in the interval $(-\frac{1}{\beta},-\alpha)$.
\end{lemma}
\begin{proof}
Let us set $z:=x+iy$ with $x,y\in\R$. Then $p(z)=0$ if and only if
  \begin{equation*}
\begin{cases}
\beta x^3+x^2+\beta b x-(3\beta x+1)y^2+a=0,\\
y(-\beta y^2+3\beta x^2+2x+\beta b)=0,
\end{cases}
\end{equation*}
from which we derive
 \begin{align}\label{sistema}
\begin{cases}
q(x):=\beta x^3+x^2+\beta b x+a=0,\\
y=0,
\end{cases}
\quad \text{or}\quad
\begin{cases}
r(x):=8\beta x^3+8x^2+2\left(\frac{1}{\beta}+\beta b\right)x+b-a=0,\\
y^2=3x^2+\frac{2}{\beta}x+b.
\end{cases}
\end{align}
By recalling $a>0$ and $b-a\geq(c_0-1)a>0$, for every $x\geq 0$ we have $q(x)>0$ and $r(x)>0$, and so the real part of the roots cannot be positive or zero. Moreover, since for every $x\leq -\frac{1}{\beta}$ we have $\beta x^3+x^2\leq 0$, we obtain
\begin{align*}
  q(x)\leq -b+a\leq (1-c_0)a<0\qquad\text{and}\qquad r(x)\leq b-a-2\left(\tfrac{1}{\beta^2}+b\right)=-b-a-\tfrac{2}{\beta^2}<0,
\end{align*}
which imply that the real part of the roots does not belong to $(-\infty,-\frac{1}{\beta}]$. Therefore, by calling $z_1,z_2,z_3\in\C$ the three roots of the polynomial $p$, we can say
\begin{equation}\label{bound-roots}
    \Re(z_i)\in(-\tfrac{1}{\beta},0)\quad \text{for $i=1,2,3.$}
\end{equation}

{\it Case 1: there is only one real root.} In this case by \eqref{sistema} there exists a unique $x_1\in (-\frac{1}{\beta},0)$ which satisfies $r(x_1)=0$ and $3x^2_1+\frac{2}{\beta}x_1+b>0$. Indeed by setting $y_1:=\sqrt{3x^2_1+\frac{2}{\beta}x_1+b}$ we obtain that $x_1+iy_1$ and $x_1-iy_1$ are two distinct non-real roots of $p$. Since 
\begin{align*}
    r(-\tfrac{1}{2\beta})&=-\tfrac{1}{\beta^2}+\tfrac{2}{\beta^2}-\tfrac{1}{\beta^2}-b+b-a=-a<0,\\
    r(-\tfrac{\beta(b-a)}{2\left(b\beta^2+1 \right)})&=\tfrac{\beta^2(b-a)^2((a+b)\beta^2+2)}{(b\beta^2+1)^3}>0,
\end{align*}
then $x_1\in(-\frac{1 }{2\beta},-\frac{\beta(b-a)}{2\left(b\beta^2+1\right)})$. Moreover
\begin{align*}
    q(-\tfrac{1}{\beta})&=-\tfrac{1}{\beta^2}+\tfrac{1}{\beta^2}-b+a=-b+a<0,\\
        q(-\tfrac{a}{\beta b})&=-\tfrac{a^3}{b^3\beta^2}+\tfrac{a^2}{b^2\beta^2}-a+a=\tfrac{a^2(b-a)}{b^3\beta^2}>0,
\end{align*}
hence there exists  $x_0\in(-\frac{1}{\beta},-\frac{a}{\beta b})$ such that $q(x_0)=0$. As a consequence of this, $(x_0,0)$ satisfies the first system in \eqref{sistema}, which implies that $x_0$ is the real root of $p$, hence we have $\Re(z_i)\in (-\frac{1}{\beta},\max\{-\frac{a}{\beta b},-\frac{\beta(b-a)}{2\left(b \beta^2+1\right)}\})$. Thanks to \eqref{ip1} we can say $-\tfrac{a}{\beta b}\leq -\tfrac{1}{c_1\beta}$ and $-\tfrac{\beta(b-a)}{2\left(b\beta^2+1\right)}\leq \tfrac{\beta(1-c_0)a}{2(c_1a\beta^2+1)}\leq\tfrac{\beta(1-c_0)a_0}{2\left( c_1a_0\beta^2+1\right)}$, where in the last inequality we use the decreasing property of the function $a\mapsto \tfrac{\beta(1-c_0)a}{2(c_1a\beta^2+1)}$. This implies
\begin{equation}\label{bound-zi-2}
\Re(z_i)\in (-\tfrac{1}{\beta},\max\{-\tfrac{1}{c_1\beta },\tfrac{\beta(1-c_0)a_0}{2\left( c_1a_0\beta^2+1\right)}\})\quad \text{for $i=1,2,3$}.    
\end{equation}

{\it Case 2: there are only real roots.} In this case we have $b\leq \frac{1}{3\beta^2}$, otherwise $q'(x)>0$ for every $x\in \R$, which forces $p$ to have also non-real roots. Thanks to \eqref{ip1} we have also $a<b\leq\frac{1}{3\beta^2}$. By setting $\tilde{b}_0:=1-\sqrt{1-3b_0\beta^2}$, we can write
\begin{align*}
    -\tilde{b}_0a_0\beta\geq-\tilde{b}_0a\beta\geq-(1-\sqrt{1-3b\beta^2})a\beta>\tfrac{-1+\sqrt{1-3 b\beta^2}}{3\beta}>-\tfrac{1}{\beta},
\end{align*}
which implies
\begin{equation}\label{q'}
     q'(x)>0\quad \text{ for every $x\in[-\tilde{b}_0a_0\beta,+\infty)$}.
\end{equation}
Since
\begin{align*}
     q(-\tilde{b}_0a_0\beta)\geq\beta^2\tilde{b}_0^2a_0^2(1-\beta^2\tilde{b}_0a_0)+a_0(1-\beta^2\tilde{b}_0b)>a_0(1+\beta^2\tilde{b}_0^2 a_0)(1-\beta^2\tilde{b}_0b)>0,
\end{align*}
thanks to \eqref{ip1}, \eqref{bound-roots}, and \eqref{q'} we get
\begin{equation}\label{bound-zi-3}
 \Re(z_i)\in (-\tfrac{1}{\beta},-\tilde{b}_0a_0\beta),\quad \text{for $i=1,2,3$}.   
\end{equation}
By combining \eqref{bound-zi-2} and \eqref{bound-zi-3}, we obtain the conclusion with $\alpha:=\min\{\tilde{b}_0a_0\beta,\frac{1}{c_1\beta },\frac{\beta(c_0-1)a_0}{2\left( c_1a_0\beta^2+1\right)}\}$.
\end{proof}

The following easy estimate is used in the proof of Lemma \ref{lem:ex-lax}.

\begin{lemma}\label{comp-in}
For every $z,w\in\C$ with $\Re(z)>0$ and $\Re(w)<0$ the following inequality holds:
\begin{equation*}
    |(z-w)(z-\bar w)|\geq |\Re(w)||\Im(w)|.
\end{equation*}
\end{lemma}
\begin{proof}
Without loss of generality we can suppose $\Im(w)>0$, otherwise we exchange the role of $w$ with $\bar w$. If $\Im(z)>0$, then
\begin{align*}
    &|z-w|\geq |\Re(z-w)|=|\Re(z)+\Re(-w)|=\Re(z)+\Re(-w)\geq|\Re(w)|, \\
    &|z-\bar w|\geq |\Im(z-\bar w)|=|\Im(z)+\Im(w)|=\Im(z)+\Im(w)\geq|\Im(w)|,
\end{align*}
which give the conclusion in this case. If $\Im(z)<0$, then
\begin{align*}
    &|z-w|\geq |\Im(z-w)|=|-\Im(-z)-\Im(w)|=\Im(-z)+\Im(w)\geq|\Im(w)|, \\
    &|z-\bar w|\geq |\Re(z-\bar w)|=|\Re(z)+\Re(-w)|=\Re(z)+\Re(-w)\geq|\Re(w)|,
\end{align*}
which conclude the proof.
\end{proof}

{\begin{acknowledgements}
This paper is based on work supported by the National Research Project (PRIN  2017) 
 "Variational Methods for Stationary and Evolution Problems with Singularities and 
 Interfaces", funded by the Italian Ministry of University and Research. 
The authors are members of the {\em Gruppo Nazionale per l'Analisi Ma\-te\-ma\-ti\-ca, la Probabilit\`a e le loro Applicazioni} (GNAMPA) of the {\em Istituto Nazionale di Alta Matematica} (INdAM).
\end{acknowledgements}}


\begin{thebibliography}{99}

\bibitem{LT} {\sc W. Arendt, C. Batty, M. Hieber, F. Neubrander}: Vector‐valued Laplace Transforms and Cauchy Problems, Monographs in Mathematics {\bf 96}. Birkhäuser, Basel, 2001.

\bibitem{DA}{\sc C. Dafermos}: An abstract Volterra equation with applications to linear viscoelasticity, {\it Journal of Differential Equations} {\bf 7} (1970), 554--569.

\bibitem{DL}{\sc R. Dautray and J.L. Lions}: Mathematical Analysis and Numerical Methods for Science and Technology, Volume 1 Physical Origins and Classical Methods, Original French edition published by Masson, S.A., Paris, 1984

\bibitem{Fab}{\sc M. Fabrizio, C. Giorgi, V. Pata}: A New Approach to Equations with Memory, {\it Arch. Rational Mech. Anal.} {\bf 198} (2010), 189-232.

\bibitem{MF-AM}{\sc M. Fabrizio and A. Morro}: Mathematical Problems in Linear Viscoelasticity. SIAM Studies in Applied Mathematics {\bf 12} 2001.


\bibitem{LM}{\sc J.L. Lions and E. Magenes}: Non-Homogeneous Boundary-Value Problems and Applications. Vol. 1 {\bf 181}, Springer-Verlag, 1972.

\bibitem{OSY}{\sc O.A. Oleinik, A.S. Shamaev, and G.A. Yosifian}: Mathematical problems in elasticity and homogenization. Studies in Mathematics and its Applications, 26. North-Holland Publishing Co., Amsterdam, 1992.

\bibitem{S}{\sc F. Sapio}: A dynamic model for viscoelasticity in domains with time--dependent cracks. Preprint SISSA 14/2020/MATE.

\bibitem{SL}{\sc L.I. Slepyan}: Models and phenomena in fracture mechanics, Foundations of Engineering Mechanics.  Springer-Verlag, Berlin, 2002.

\end{thebibliography}
\end{document}